\documentclass[11pt,a4paper]{article}
\usepackage{amsmath,amsfonts,amssymb,mathrsfs}
\usepackage{mathtools}
\usepackage{graphicx}
\usepackage{xfrac}
\usepackage{subfigure}
\usepackage{pstricks}
\usepackage{pst-node}
\usepackage{xcolor}
\usepackage[english]{babel}
\usepackage{mathscinet}

\newtheorem{theorem}{\bf Theorem}[section]
\newtheorem{lemma}[theorem]{\bf Lemma}
\newtheorem{corollary}[theorem]{\bf Corollary}
\newtheorem{definition}[theorem]{\bf Definition}
\newtheorem{remark}[theorem]{\bf Remark}
\newtheorem{proposition}[theorem]{\bf Proposition}

\newenvironment{proof}{\noindent {\sc Proof.}}{\hfill $\square$}

\def \GK {{\Gamma_{\!\! K}}}
\def \GL {{\Gamma_{\!\! L}}}
\newcommand{\R}{\mathbb{R}}

\newcommand{\N}{\mathbb{N}}
\newcommand{\E}{\mathbb{E}}
\newcommand{\QQ}{\mathbb{Q}}

\newcommand{\I}{\mathbb{I}}
\newcommand{\GG}{\mathbb{G}}
\newcommand{\KK}{\mathbb{K}}

\def \L {\mathscr{L}}
\def \K {\mathscr{K}}

\def \a {{\alpha}}
\def \b {{\beta}}

\def \g {{\gamma}}
\def \d {{\delta}}
\def \e {{\varepsilon}}

\def \epsilon {{\varepsilon}}

\def \l {{\lambda}}
\def \r {{\rho}}
\def \s {{\sigma}}
\def \t {{\tau}}
\def \m {{\mu}}

\def \x {{\xi}}
\def \y {{\eta}}
\def \z {{\zeta}}
\def \w {{\omega}}
\def \phi {{\varphi}}

\def \G {{\Gamma}}
\def \O {{\Omega}}

\def \loc {{\text{\rm loc}}}

\def\p{\partial}

\def \Q {{Q}}

\usepackage{mathtools}
\DeclarePairedDelimiter{\abs}{\lvert}{\rvert}

\begin{document}
	\title{Existence of a Fundamental Solution of Partial Differential Equations associated to Asian Options}
	\author{{\sc{Francesca Anceschi}
	\thanks{Dipartimento di Scienze Fisiche, Informatiche e Matematiche, Universit\`{a} 
	degli Studi di Modena e Reggio Emilia, via
	Campi 213/b, 41125 Modena (Italy). E-mail: francesca.anceschi@unimore.it}
	{\sc{Silvia Muzzioli}
	\thanks{Dipartimento di Economia ``Marco Biagi", Universit\`{a} degli Studi di Modena e 
	Reggio Emilia, viale Berengario 51, 41121 Modena (Italy). E-mail: 
	silvia.muzzioli@unimore.it}
	\sc{Sergio Polidoro}
	\thanks{Dipartimento di Scienze Fisiche, Informatiche e Matematiche, Universit\`{a} 
	degli Studi di Modena e Reggio Emilia, Via Campi 213/b, 41125 Modena (Italy). E-mail: 
	sergio.polidoro@unimore.it}
	}}}
	
	\date{ }
	
	\maketitle

	\bigskip

	\begin{abstract}
	\noindent
		We prove the existence and uniqueness of the fundamental solution for Kolmogorov operators associated to 
		some stochastic processes, that arise in the Black \& Scholes setting for the pricing problem 
		relevant to path dependent options. We improve previous results in that we provide a closed form expression
		for the solution of the Cauchy problem under weak regularity assumptions on the coefficients of the differential operator.
		Our method is based on a limiting procedure, whose convergence relies on some barrier arguments and uniform a priori estimates recently discovered.
	
		
	\end{abstract}
	
\setcounter{equation}{0}
\section{Introduction}
Asian options belong to the family of \emph{path-dependent options} whose payoff depends on the average of the underlying stock 
price over a certain time interval. In the Black \& Scholes framework, the price of the underlying Stock $S_t$ and of the bond $B_t$ are described by the processes
\begin{equation*}
 S_t = S_0 e^{\mu t + \sigma W_t}, \qquad B_t = B_0e^{rt}, \qquad 0 \le t \le T, 
\end{equation*}
where $\mu, r, T$, and $\sigma$ are given constants. If the price observations are considered as a set of regularly spaced time points we refer to a discrete Asian Option, otherwise when we consider a continuum of price observations and its average it is computed by means of an integral we have a continuous Asian Option. In particular, in this work we consider continuous Asian Options. In the Black \& Scholes setting, the price $(Z_t)_{0 \le t \le T}$ of a path dependent option is considered as a function $Z_{t} = Z(S_{t},A_{t},t)$ that depends on the \emph{stock price} $S_t$, the \emph{time to maturity} $t$ and on an average$A_t$ of the stock price
\begin{equation} \label{At}
   A_{t} \, = \, \int \limits_{0}^{t} f(S_{\t}) \, d\t, \qquad t \in [0, T].
\end{equation}
From a financial point of view, Asian Options have several advantages. Indeed, they are less expensive than Plain Vanilla 
Options thanks to the averaging mechanism which allows to reduce their volatility as well. Secondly, they reduce the risk of 
market manipulation of the underlying instrument at maturity (see \cite{sun2015asian}). In this sense, Asian Options are suitable 
to fulfill some of the needs of corporate treasures.

We refer to the Black \& Scholes \cite{BS} and to Merton \cite{ME} articles for the seminal works of this theory, and to the books by Bj\"ork \cite{B}, Hull \cite{HU} and Pascucci \cite{PA} for a comprehensive treatment of the recent development of this subject. The most common techniques to price path-dependent derivatives are:
\begin{itemize}
\item \textbf{the Monte Carlo simulations}, relying on the Feynman–Kac formula
\begin{equation}
	Z(S,A,t) \, = \, \E^{\QQ} \left[ e^{ - r(T - t)} \phi(S_{T} ,A_{T} ) \, | \, 
	(S_{t},A_{t}) \, = \, (S,A) \right] ,
\end{equation}
where $\QQ$ is a measure such that the process $e^{- rt} Z_{t}$ is a martingale under 
$\QQ$, and the fast Fourier transform (see for example \cite{GY}, \cite{ballestra2007numerical}, 
\cite{fu1999pricing}). In \cite{GY}, the authors derived an analytical expression for the Laplace transform in maturity for 
the continuous call option case when the asset price follows a geometric Brownian motion. However, as pointed out by \cite{foschi2013approximations}, \cite{Dufresne} and \cite{fu1999pricing} the analytical method of \cite{GY} 
can lead to numerical problems for short maturities or small volatilities. These problems are consequences of the slowly 
decaying oscillatory nature of the integrand for such parameter values (see \cite{fu1999pricing}).
\item \textbf{The PDE approach}, which has the aim to solve numerically the Cauchy problem associated with the no-arbitrage PDE. 
Related works following this line are those of \cite{CPR}, \cite{BPV}, \cite{vecer2001new}. In \cite{curran1994valuing}, the author applies a method on conditioning on the geometric mean price. In this case an approximation of Arithmetic Asian option prices is available. In \cite{dewynne2008differential}, the author derives an accurate approximation formulae for Asian-rate Call options in the Black \& Scholes model by a matched asymptotic expansion. In this work we rely on the results proved in \cite{CPR}, where the authors prove via probabilistic techniques the existence of the fundamental solution $\G$ for the operator $\L$ with smooth coefficients $a$ and $b$. 
Moreover, we recall the existence and local regularity results proved by Lanconelli, Pascucci and Polidoro \cite{LPP}, 
under the assumption that the coefficients $a$ and $b$ belong to some space of H\"older continuous functions. 
\end{itemize}
In this work we consider the analytical approach based on PDEs since it has several advantages compared to the Monte Carlo approach, as it is stressed by \cite{foschi2013approximations}. First of all, providing an analytical approximation of the solution in closed-form gives evidence of the explicit dependency of the results on the underlying parameters. Secondly, analytical approaches produce better and faster sensitivities than Monte Carlo methods. 

In order to explain our main results, we introduce some notation. From now on, we consider the stochastic differential equation of the process $(S_{t},B_{t},A_{t})_{t \ge 0}$
\begin{equation} \label{SBA-proc}
   \begin{cases}
   	dS_{t} \, = \, \m(S_{t}, A_{t}, t) \, S_{t} dt \, + \, \s(S_{t}, A_{t}, t) \, S_{t} dW_{t}, \\
	dB_{t} \, = \, r(S_{t}, A_{t}, t) \, B_{t} \, dt \\
	dA_{t} \, = \, f (S_{t}) \, dt, 
   \end{cases}
\end{equation}
where $t \in ]0,T[$, $\m$, $r$ and $\s$ depend on $S_{t}$, $A_{t}$ and $t$. Then we construct the replicating portfolio $(Z_t)_{0 \le t \le T}$ for the option, we consider it as a function $Z_{t} = Z(S_{t},A_{t},t)$ and we apply It\^o’s formula. Thus, we obtain the following Cauchy problem
\begin{equation} \label{PDE1}
	\begin{cases}
		\tfrac12 \s^{2}(S, A, t) S^{2} \frac{\p^{2} Z}{\p S^{2}} + f(S) \frac{\p Z}{\p A} + r(S, A, t) \left( S \frac{\p Z }{\p S} - Z \right) +
		 \frac{\p Z}{\p t} = 0 \quad & (S, A, t) \in \R^{+} \times \R^{+} \times ]0, T[, \\
		Z(S,A,T) = \phi (S, A)   &(S, A) \in \R^{+} \times \R^{+},
	\end{cases}
\end{equation}
where $\phi$ is the \emph{payoff} of the Asian Option. 

We emphasize that the PDE approach adopted in this work improves the previously known results in that it allows us to consider differential operators with \emph{locally H\"older continuous} coefficients, which is a milder assumption than the usual ones. We will be more specific in the following, as we introduce the required notation. Our approach can be also applied to more general problems than the one described above. For instance, in a further investigation we will consider the pricing problem for an Option on a basket containing $n$ assets $S_{t} = (S_{t}^{1}, \ldots, S_{t}^{n})$ whose dynamic is
\begin{equation}
	dS_{t}^{j} \, = \, S_{t}^{j} \, \m_{j}(S_{t},A_{t},t) \, + \, S_{t}^{j} \, 
	\sum \limits_{k=1}^{n} \s_{jk} (S_{t}, A_{t}, t ) \, d W_{t}^{k}, 
	\qquad j = 1, \ldots, n, 
\end{equation}
where $(W_{t}^{1},..., W_{t}^{n})_{t \ge 0}$ is a $n$-dimensional Wiener process and $(A_{t})_{t\ge0}$ is an average of the assets. 

%

\subsection{Geometric Average Asian Options}
We deal with Geometric Average Asian Options when we choose $f(S) = \log(S)$ in the formula \eqref{At}. Through a simple change of variable $v (x, y, t) := Z(e^{x}, y, T - t)$ we transform the PDE \eqref{PDE1}, with its final condition, into the following Cauchy problem
\begin{equation} \label{PDE2}
	\begin{cases}
		\tfrac12 \s^{2}(x,y,t) \left( \frac{\p^{2}v}{\p x^{2}} - \frac{\p v}{\p x} \right) 
		\, + \, x \frac{\p v}{\p y} \, + \, r(x,y,t) \left( \frac{\p v}{\p x} - v \right) \, = \, \frac{\p v}{\p t} \\
		v(x,y,0) = \widetilde \phi (x,y),
	\end{cases}
\end{equation}
where $\widetilde \phi (x,y) := \phi(e^{x},y)$.
Note that, if we assume that $\frac{\partial \s}{\partial x}$ is a continuous function, then the differential operator in \eqref{PDE2} can be written in its \emph{divergence form}. Precisely, for every sufficiently smooth function $u$, we have that the PDE in \eqref{PDE2} writes as $\K u = 0$, with
\begin{equation} \label{K}
	\K u (x,y,t) \, = \, \frac{\p}{\p x} \biggl(a (x,y,t) \frac{\p u}{\p x} \biggr) \, + \, b(x,y,t) \, \frac{\p u}{\p x} \, + \,  x \frac{\p u}{\p y} \, - r(x,y,t) u - \, \frac{\p u}{\p t}.
\end{equation}
Here $a(x,y,t) = \tfrac12 \s^{2}(x,y,t)$ and $b(x,y,t) = r(x,y,t) - \tfrac12 \s^{2}(x,y,t) - \s(x,y,t)\frac{\p \s(x,y,t)}{\p x}$. The reason to write $\K$ in this form is that we need apply some results that have been proved only for divergence form operators. 
We also introduce its formal adjoint $\K^{*}$, acting on smooth functions $w=w(x,y,t)$ as 
\begin{equation}
	\label{Kadj}
	\K^{*} w(x,y,t) \, = \, \frac{\p}{\p x} \biggl(a (x,y,t) \frac{\p w}{\p x} \biggr) \, - \, \frac{\p}{\p x} \biggl( b(x,y,t) \, w \biggr) \, - \, x \frac{\p w}{\p y} \, - \, r(x,y,t) w \, + \, \frac{\p w}{\p t}.
\end{equation}

The Cauchy problem \eqref{PDE2} has been studied over the years, and the fundamental solution $\GK$ associated to the operator $\K$ provides us with a representation formula for its solution. In particular, if $\widetilde \phi$ is a bounded continuous function, then 
\begin{equation} \label{rep1}
	u(x,y,t) \, = \, \int \limits_{\R^{+} \times \R} \GK(x,y,t; \x, \y, 0) \, \widetilde \phi(\x, \y) \, d\x \, d\y, 
\end{equation}
is a classical solution to \eqref{PDE2}. Kolmogorov wrote in \cite{K1} the explicit expression of the fundamental solution $\GK$ for the operator $\K$ with constant coefficients $\sigma$ and $r$. In this case the function $u(x,y,t) := e^{rt} v \left( x + \left( \tfrac12 \s^{2} - r \right) t,y,t \right)$ is a solution to the Cauchy problem 
\begin{equation*}
 	\begin{cases}
 		\K_\l u = 0, \\
 		u(x,y,0) = \widetilde \phi \left(x,y\right),
 	\end{cases}
\end{equation*}
where $\l = \tfrac12 \s^{2}$ and 
\begin{equation} \label{Klambda}
	\K_\l := \l \, \frac{\p^{2}}{\p x^{2}} + x \frac{\p }{\p y} - \frac{\p u}{\p t},
\end{equation}
Moreover, the fundamental solution $\GK^{\!\!\!\! \l}$ of the operator $\K_{\l}$ is
\begin{equation} \label{Glambda}
	\GK^{\!\!\!\! \l}(x,y,t,\x,\y,\t) \, = \,
	\begin{cases}
	\frac{\sqrt3}{2 \l \pi (t - \t)^{2}} \exp \left( - \frac{\big \lvert x - \x \big \rvert^{2}}{4 \l (t-\t)} 
		- 3 \frac{\big \lvert y - \y + (t-\t) \tfrac{x + \x}{2} \big \rvert^{2} }{ \l (t- \t)^{3}} \right) \quad &t > \t, \\
	0 \qquad &t \le \t.
	\end{cases}
\end{equation}
Thus, we obtain a closed form for the price of the Geometric Average Asian Option in the case of constant volatility $\sigma$ and interest rate $r$.


The Levy parametrix method provied us with a fundamental solution for operators in the form $\K$ with H\"older continuous coefficients. This method has been used by several authors \cite{P2, PDF, PODF, KMM} and requires a 
uniform H\"older continuity of the coefficients of $\K$. The definition of the H\"older space $C^{\a}_K$ is given in Section \ref{holder}. We will see that a function $f$ belongs to the space $C^{\a}_K$, with $0 < \alpha \le 1$, if there exists a positive constant $M$ such that
\begin{equation} \label{calpha}
	|f(x,y,t) - f(\x, \y, \t))| \, \le M \left( \abs{x - \x} \, + \, \Big \lvert y - \y + (t- \t)\tfrac{x + \x}{2} \Big \rvert^{\tfrac13} \, + \, \abs{t-\t}^{\tfrac12}\right)^\alpha,
\end{equation}
for every $(x,y,t), (\x, \y, \t) \in \R^3$.
On one hand, the intrinsic H\"older space $C^{\a}_K$ associated to $\K$ complies with the fundamental solution 
$\G^{\l}$ written in \eqref{Glambda}. 
On the other hand, intrinsic H\"older regularity can be a rather restrictive property, as it has already been pointed out by Pascucci and Pesce in the Example 1.3 of \cite{PaPe}. In particular, Pascucci and Pesce show that, whenever a function $f = f(y)$ only depends on $y$ and belongs to $C^{\a}_K$, is necessarily constant. As we will see in the sequel of this article, we only require a \emph{local} H\"older regularity of the coefficients of the operator $\K$. This allows us to consider a wider family of continuous functions. More precisely, we consider the 
following assumption on the coefficients $a$ and $b$:
\begin{itemize}
 	\item[\textbf{(H$_{K}$)}] $a, b, r$, $\frac{\partial a}{\partial x}, \frac{\partial b}{\partial x}$ $\in C^{\a}_{\loc}(\R^{3})$. Moreover, there exist two positive constants $\l$, $\Lambda$ such that 
		\begin{align*}
			\l \le a(x,y,t) \le \Lambda, \quad |b(x,y,t)|, |r(x,y,t)|, \left| \tfrac{\partial a}{\partial x}(x,y,t) \right|, \left| \tfrac{\partial b}{\partial x}(x,y,t) \right| \le \Lambda, 
		\end{align*}
for every $(x,y,t) \in \R^{3}$.
\end{itemize}
In the above display, $C^{\a}_{\loc}$ denotes the usual space of H\"older continuous functions. Proposition \ref{manfredini} states that a function $f$ belongs to the space of \emph{locally H\"older continuous} functions $C^{\a}_{K, \loc}$ if, and only if, it belongs to the space $C^{\b}_{\loc}$ for some positive $\b$.
We are now in position to state our main results. 
\begin{theorem}
	\label{main1}
Let us consider the operator $\K$ under the assumption \textbf{(H$_{K}$)}. Then there exists a unique fundamental solution $\GK$ of $\K$ in the sense of Definition \ref{funsolk}.
%
Moreover, the following properties hold:
	\begin{enumerate}
		\item Support of $\GK$: for every $(x,y,t)$, $(\x, \eta, \t) \in \R^{3}$ with $t \le \t$
			\begin{equation*}
				\GK (x,y,t; \x, \y, \t) = 0;
			\end{equation*}
		\item Reproduction property: for every $(x,y,t), (x_{0}, y_{0}, t_{0}) \in \R^{3}$ and $\tau \in \R$ with $t_{0} < \t < t$
			\begin{equation*}
				\GK (x,y,t; x_{0}, y_{0}, t_{0}) = \int \limits_{\R^{2}}
				\GK (x,y,t; \x, \y, \t) \, \GK(\x, \y, \t; x_{0}, y_{0}, t_{0}) \, d\x \, d\y;
			\end{equation*}
	\item Integral of $\GK$: the following bound holds true
		\begin{equation*}
             e^{-\Lambda (t - \tau)} \, \le \, 
 			\int \limits_{\R^{2}} \GK(x,y,t; \x, \y, \t) \; d \x \, d\y 
 			\, \le \, e^{\Lambda (t - \tau)};
 		\end{equation*}
	\item Bounds for $\GK$: let $I=]T_{0}, T_{1}[$ be a bounded interval, then there exist four positive constants 
		$\l^{+}$, $\l^{-}$, $C^{+}$, $C^{-}$ such that for every $(x,y,t), (\x, \y, \t) \in \R^{3}$ with $T_{0} < t < T < T_{1}$
       \begin{align*}
		C^{-} \, \GK^{\! \! \! \! \l^{-}} (x,y,t; \x, \y, \t) \, \le \, \GK(x,y,t; \x, \y, \t) \, \le \,
		C^{+} \, \GK^{\! \! \! \! \l^{+}} (x,y,t; \x, \y, \t).
		\end{align*}
		The constants $\l^{+}$, $\l^{-}$, $C^{+}$, $C^{-}$ depend only on $\K$ and $T_{1}-T_{0}$.
		$\GK^{\! \! \! \! \l^{-}}$ and $\GK^{\! \! \! \! \l^{+}}$ respectively denote the fundamental solution of 
		$\K_{\l^{-}}$ and $\K_{\l^{+}}$, defined in \eqref{Glambda} and \eqref{Klambda} respectively.
	\end{enumerate}
	Moreover, the function $\GK^{*}(\x, \y, \t; x, y, t) = \GK(x,y,t; \x, \y, \t)$ is the fundamental solution of the adjoint 
	operator $\K^{*}$ with pole at $(\x, \y, \t)$ and satisfies all of the previous properties accordingly.
\end{theorem}

We remark that in the Black\&Scholes setting it is natural to consider the Cauchy problem \eqref{PDE1} with an unbounded initial condition $\varphi$ that grows linearly. 
After the change of variable $v (x, y, t) := Z(e^{x}, y, T - t)$, it corresponds to an \emph{exponential growth} for the Cauchy data $\widetilde \varphi$. As we will see in Remark \ref{rem-growht}, the formula \eqref{rep1} supports initial data satisfying the following condition
\begin{equation} \label{eq-growht}
 | \widetilde \varphi (x,y) | \le M \exp(C |(x,y)|^\alpha), \qquad (x,y) \in \R^2,
\end{equation}
for some positive constants $M, C$ and $\alpha$, with $\alpha < 2$. 
Note that if we consider $\a=2$ then the solution to the Cauchy problem \ref{PDE1} is defined in a suitably small interval of time. 
Moreover, the following uniqueness result holds. 

\begin{theorem} \label{uniquenessK}
Let us consider the operator $\K$ under the assumption \textbf{(H$_{K}$)}. 
Let $u_{1}$ and $u_{2}$ be classical solutions to 
 	\begin{equation} \label{cauchyproblem}
		\begin{cases}
    				\K u = 0, \qquad & (x,y,t) \in \R^{2} \times ]t_0, T] \\
    				u(x_0,y_0,t_0)=\varphi (x_0,y_0) &(x_0,y_0) \in \mathbb{R}^2.
 			\end{cases}
	\end{equation} 
in the sense of Definition \ref{solution}, and
\begin{equation*}
	|u_{1}(x,y,t)| + |u_{2}(x,y,t)| \, \le \,M \exp\left(C (x^{2} + y^{2})\right),
\end{equation*}
for some positive consants $M$ and $C$, then $u_{1} = u_{2}$ in $\R^{2} \times ]t_{0},T]$.
\end{theorem}

\subsection{Arithmetic Average Asian Options}
When we deal with $f(S) = S$, we consider Arithmetic Average Asian Options. Through the change of variable $v (x, y, T - t) := Z(x, y, t)$ we transform the Cauchy problem \eqref{PDE1} into the following
\begin{equation*}
	\begin{cases}
		\tfrac12 \s^{2}(x,y,t)  x^{2} \frac{\p^{2} v}{\p x^{2}} \,  
		+ \, r(x,y,t) \left( x \frac{\p v}{\p x} \, - \, v \right) \, + \, x \frac{\p v}{\p y} (x,y,t) \, = \, \frac{\p v}{\p t} \\
		v(x,y,0) = \phi (x, y)   .
	\end{cases}
\end{equation*}
As we did for the case of Geometric Average Asian Options, we write the operator appearing in the above PDE in its divergence form 
\begin{equation}
	\label{PDE3}
	\begin{cases}
		x \, \frac{\p}{\p x} ( a (x,y,t) \, x \, \frac{\p v}{\p x} ) \, + \, b(x,y,t) \, x  \, \frac{\p v}{\p x}  \, 
		- r(x,y,t) v \, + \, x \, \frac{\p v}{\p y} = \, \frac{\p v}{\p t}, \\
		v(x,y,0) = \phi (x, y),
	\end{cases}
\end{equation}
where $a(x,y,t) = \tfrac12 \s^{2}(x,y,t)$ and $b(x,y,t) = r(x,y,t) - 2 x a(x,y,t) - x^2 \p_x a(x,y,t)$. Note that the coefficient $x$ in front of the derivative $\tfrac{\partial}{\partial x}$ introduces new difficulties. We denote by $\L$ the differential operator in \eqref{PDE3}, acting on sufficiently smooth functions $u = u(x,y,t)$ as 
\begin{equation} \label{Ldiv}
	\L u(x,y,t) \, := \, x \, \frac{\p }{\p x} \left( a (x,y,t) \, x \, \frac{\p u}{\p x} \right) \, + \, b(x,y,t) \, x \, \frac{\p u}{\p x}  \, + \, x \, \frac{\p u}{\p y} - r(x,y,t) u - \, \frac{\p u}{\p t},
\end{equation}
and its formal adjoint $\L^*$, acting on differentiable functions $w=w(x,y,t)$ as follows
\begin{equation} \label{L*}
	\L^* w (x,y,t) \, := \, x \, \frac{\p w}{\p x} \left( a (x,y,t) \, x \, \frac{\p w}{\p x} \right) \, 
	- \, \frac{\p}{\p x}  \biggl( x \, b(x,y,t) \, w  \biggr)  \, - \, x \, 
	\frac{\p}{\p y} w \, - r(x,y,t) w + \, \frac{\p}{\p t} w.
\end{equation}
The simplest form of the operator $\L$ is associated to the stochastic process 
$(S_t , A_t )_{t\ge 0}$
\begin{equation} \label{processY}
	S_{t}=S_0 e^{\mu t + \sigma W_t}, \qquad A_t = \int_0^t \exp( \m s + \s W_s) ds,
\end{equation}
where $W_{t}$ is a real valued Brownian motion starting from $0$. Indeed, when $\m = 0$ and $\s$ 
is a positive constant we can consider the following model operator
\begin{equation} \label{Llambda}
	\L_\l \, :=  \l \, x^{2} \p^{2}_{x} \, + \, x \p_{x} \, + \, x \p_{y} \, - \, \p_{t},
\end{equation}
where $\l=\tfrac12 \s^{2}$.
As it is pointed out by Yor in \cite{Y}, thanks to the scaling invariance properties of the Brownian motion we can restrict 
ourselves to the case where $\s = \sqrt{2}$, for which he proves the existence of the transition density of 
the associated process $(S_t , A_t )_{t\ge 0}$, which reads as follows
\begin{equation*}
	p(w,y,t) = \frac{e^{\tfrac{\pi^{2}}{2t}}}{\pi \sqrt{2\pi t}} \, \exp \left( - \frac{1+e^{2w}}{2y} \right)
	\, \frac{e^{w}}{y^{2}} \, \Theta \left( \frac{e^{w}}{y}, t \right),
\end{equation*}
where
\begin{equation*}
	\Theta (z,t) = \int \limits_{0}^{\infty} e^{- \tfrac{\x^{2}}{2t}} e^{-z \cosh(\x)} \sinh (\x) \, 
			\sin \left( \frac{\pi \x}{t} \right) \, d\x.
\end{equation*}
Thus, the explicit expression of the fundamental solution $\GL^{\!\! 1}$ associated with the operator
\begin{equation} \label{L0}
	\L_1 \, := \, x^{2} \p^{2}_{x} \, + \, x \p_{x} \, + \, x \p_{y} \, - \, \p_{t},
\end{equation}
reads as follows
\begin{equation}
	\label{fun-yor}
	\GL^{\!\! 1} (x,y,t;x_{0}, y_{0}, t_{0}) = \frac{1}{2x x_{0}} \, p \left(\frac12 \, \log \left( \frac{x_{0}}{x} \right), 
	\frac{y-y_{0}}{2x}, \frac{t-t_{0}}{2} \right).
\end{equation}
Thus, the pricing problem for the simplest case of Arithmetic Average Asian Options can be solved by the argument outlined in the previous subsection when considering \eqref{PDE2}, but in this case the differential operator $\K$ needs to be replaced by $\L_{\l}$. 
As we can see from the explicit expression \eqref{fun-yor} of the fundamental solution $\GL^{\!\!\!\! \l}$ of $\L_{\l}$, and as several authors point out (for instance, see \cite{AimiDiazziGuardasoni, Dufresne, FoschiPagliaraniPascucci, FuMadanWang, Shaw}), the explicit representation of the Asian option prices given by Geman and Yor in \cite{GemanYor} is hardly numerically treatable, in particular when pricing Asian Options with short maturities or small volatilities.
We quote \cite{Yor2, GemanYor} for an exhaustive presentation of the topic, other related works are due to Matsumoto, Geman and Yor \cite{Matsumoto, GemanYor, MatsumotoYor}, Carr and Schr\"oder \cite{CarrScroder}, Bally and Kohatsu-Higa \cite{BallyKohatsuhiga}. 


As we have already pointed out at the beginning of the Introduction, in this work we consider the operator $\L$ with variable coefficients. This allows one to deal with more general market models, but the mathematical theory for this kind of operator $\L$ is nowadays still incomplete. 
Indeed, the unique result available on the existence of a fundamental solution for $\L$ has been proved by Cibelli, Polidoro and Rossi in \cite{CPR} and requires the $C^\infty$ smoothness of the coefficients $a$ and $b$. Moreover, only the case $r=0$ is considered in \cite{CPR}.
Our research weakens the regularity requirements on the coefficients in that only the \textit{local} H\"older continuity is needed to produce classical solutions to the pricing problem. The class of H\"older continuous functions $C^{\a}_{L}$ related to the operator $\L$ is strongly linked to the definition of the space $C^{\a}_{K}$ related to the operator $\K$, as we will see in the sequel of this article. 
Moreover, in the following we prove that \textit{locally} the two definitions coincide (see Proposition \ref{coincide}).
This allows us to consider a wider family of continuous functions, since the local H\"older condition is really easy to check 
and less restrictive that the \emph{global} H\"older continuity assumption, required for instance by the parametrix method, that is an alternative approach to produce a fundamental solution. We are now ready to state the precise assumption for the coefficients $a$ and $b$ of the operator $\L$.
\begin{itemize}
 	\item[\textbf{(H$_{L}$)}] $a, b$, $\frac{\partial a}{\partial x}$, $\frac{\partial b}{\partial x}$ $\in C^{\a}_{\loc}(\R^{+} 		\times \R^2)$. Moreover, there exist two positive constants $\l$, $\Lambda$ such that 			
		\begin{align*}
			\l \le a(x,y,t) \le \Lambda, \quad  \abs{b(x,y,t)}, \left| \tfrac{\partial}{\partial x} \left( x \, a(x,y,t) 
			\right) \right|, \left| \tfrac{\partial}{\partial x}  \left( x \, b(x,y,t) \right) \right| \le \Lambda,
		\end{align*}
	for every $(x,y,t) \in \R^{+} \times \R^2$.
\end{itemize}

\begin{remark} 
 As said above, we only consider the operator $\L$ without the zero order term $r$, as we rely on the results proved in \cite{CPR}, where this condition was assumed. However a simple change of function allows us to consider any continuous function $r = r(t)$ only depending on $t$. Indeed, if $u$ is a solution to $\L u = 0$, where the term $r$ doesn't appear in $\L$, then the function
 \begin{equation*}
				v(x,y,t) = e^{R(t)} u(x,y,t), \qquad R(t) = \int_{t_0}^{t} r(s) d s, 
\end{equation*}
solves the equation $\L v (x,y,t) - r(t) v(x,y,t)= 0$. 
\end{remark}

We are now in position to state our main result. 
\begin{theorem}
	\label{main2}
	Let us consider the operator $\L$ under the assumption \textbf{(H$_{L}$)}. Then there exists 
	a unique fundamental solution $\GL$ of $\L$ in the sense of Definition \ref{funsoll}.
%
Moreover, the following properties hold:
	\begin{enumerate}
		\item Support of $\GL$: for every $(x,y,t)$, $(\x, \eta, \t) \in \R^{+} \times \R^2$ with $t \le \t$ or $y \ge \y$
			\begin{equation*}
				\GL (x,y,t; \x, \y, \t) = 0;
			\end{equation*}
		\item Reproduction property:  for every $(x,y,t)$, $(x_{0}, y_{0}, t_{0}) \in \R^{+} \times \R^2$ and $\tau \in \R$ 
		with $t_{0} < \t < t$
			\begin{equation*}
				\GL (x,y,t; x_{0}, y_{0}, t_{0}) = \int \limits_{\R^{+} \times \R}
				\GL (x,y,t; \x, \y, \t) \, \GL(\x, \y, \t; x_{0}, y_{0}, t_{0}) \, d\x \, d\y;
			\end{equation*}
		\item Integral of $\GL$: there exists a positive constant $\overline{C}$ depending on $t-\t$ and such that 
			$\overline{C} \rightarrow 1$ as $t \rightarrow \t$ for which
			\begin{equation*}
				\int \limits_{\R^{+} \times \R} \GL (x,y,t; \x, \y, \t) \, d\x \, d\y  \, = \, 1,     \qquad 
				\int \limits_{\R^{+} \times \R} \GL (x,y,t; \x, \y, \t) \, dx \, dy  \, = \, \overline{C};
			\end{equation*}
		\item Bounds for $\GL$: for every $\e \in ]0,1[$, and $T>0$ there exist two positive constants $c_{\e}^{-}, C_{\e}^{+}$ 
			depending on $\e$, on $T$ and on the operator $\L$, and two positive constants $C^{-}, c^{+}$, only depending on the 
			operator $\L$, such that 
			\begin{align*}
				\frac{c_{\e}^{-}}{x_{0}^{2} (t - t_{0})^{2}} & \exp \left( -C^{-} \psi (x,y + x_{0} \e (t-t_{0}) , t - \e (t-t_{0}); 
				x_{0}, y_{0},t_{0}) \right) \le \\ \nonumber
				& \qquad \le \GL (x,y,t; x_{0}, y_{0},t_{0}) \le \frac{C_{\e}^{+}}{x_{0}^{2} (t - t_{0})^{2}} \exp 
				\left(-c^{+} \psi (x,y - x_{0} \e, t + \e; x_{0}, y_{0},t_{0}) \right)
			\end{align*}
			for every $(x,y,t) \in \R^{+} \times ]-\infty, y_{0} - x_{0}\e(t-t_{0})[ \times ]t_{0}, T]$, where $\psi$ is the value function
			for the optimal control problem \eqref{psi}.
	\end{enumerate}
	Moreover, the function $\GL^{*}(\x, \y, \t; x, y, t) = \GL(x,y,t; \x, \y, \t)$ is the fundamental solution of the adjoint 
	operator $\L^{*}$ with pole at $(\x, \y, \t)$ and satisfies all of the previous properties accordingly.
\end{theorem}
We note that the upper and lower bounds in the above point \emph{4.} can be written in terms of the function \eqref{fun-yor} as stated in Corollary \ref{corollaryL} below. As in the case of Geometric Average Asian Options, we consider the Cauchy problem \eqref{PDE1} with an initial condition $\varphi$ that grows linearly. However, in the case of Arithmetic Average Asian Options the change of variable $v (x, y, t) = Z(e^{x}, y, T - t)$  doesn't simplify the proof of our results. Therfore we don't apply it and we keep the linear growth as the natural assumption on the function $\varphi$. We will see in Remark \ref{rem-growht-L} that the formula \eqref{rep2} supports initial data satisfying this condition. As far as we are concerned with the uniqueness of the solution to the Cauchy problem for operators of the form \eqref{Ldiv}, we have the following result.

\begin{theorem} \label{uniquenessL}
Let us consider the operator $\L$ under the assumption \textbf{(H$_{L}$)}. 
Let $u_{1}$ and $u_{2}$ be classical solutions to 
 	\begin{equation*} \label{cauchyproblemL}
		\begin{cases}
    				\L u = 0, \qquad & (x,y,t) \in \R^{+} \times \R \times ]t_0, T] \\
    				u(x_0,y_0,t_0)=\varphi (x_0,y_0) &(x_0,y_0) \in \R^{+} \times \R.
 			\end{cases}
	\end{equation*} 
in the sense of Definition \ref{solutionL}, and
\begin{equation*}
	|u_{1}(x,y,t)| + |u_{2}(x,y,t)| \, \le \,M \exp\left(C (\log(x^{2} + y^{2} + 1) - \log(x)) + 1\right)^2,
\end{equation*}
for some positive consants $M$ and $C$, then $u_{1} = u_{2}$ in $\R^{+} \times \R \times ]t_{0}, T]$.
\end{theorem}

\bigskip

This article is structured as follows. 
In Section \ref{holder} we recall some known facts about the regularity and invariance properties of the Kolmogorov operator $\K$, in Section \ref{preliminaries} we collect known facts on the operator $\L$. Section 4 is devoted to the proof of our main results. 

\setcounter{equation}{0} 
\section{Geometric setting and fundamental solution for $\K$}
\label{holder}
This section is devoted to the study of the geometric setting suitable for the study of the Kolmogorov operator $\K$, 
and to recall some well known results concerning its fundamental solution $\GK$.
Moreover, we discuss in detail the space of locally H\"older continuous functions considered 
in the assumption \textbf{(H$_{K}$)}. For a comprehensive treatment of this subject we refer to the survey paper \cite{APsurvey}.

Firstly, let us consider the operator $\K_1$ defined in \eqref{Klambda} with $\l = 1$:
\begin{equation} \label{K0}
	\K_{1} \, := \, \frac{\p^{2}}{\p x^{2}} + x \frac{\p }{\p y} - \frac{\p u}{\p t}.
\end{equation} 
Even tough it is a strongly degenerate operator, it is \textit{hypoelliptic} in the following sense. Let $u$ be a 
distributional solution of $\K_1u = f$ in $\O \subset \R^{3}$, then
\begin{equation}
	\label{hypo}
		u \in C^{\infty}(\O) \quad \text{whenever} \quad f \in C^{\infty}(\O).
\end{equation}
H\"ormander introduced in his seminal paper \cite{H} a simple sufficient condition to check the hypoellipticity of
any second order linear differential operator defined on some open set $\O \subset \R^{N+1}$ that can be written as 
a sum of squares of smooth vector fields $X_{0}, X_{1}, \ldots, X_{m}$, as follows
\begin{equation}
	\label{sumsquares}
	\sum \limits_{i=1}^{m} X_{i}^{2} \, + \, X_{0}.
\end{equation}
 
The celebrated hypoellipticity result due to H\"ormander reads as follows. 

\medskip

\noindent
\textsc{Theorem (H\"ormander hypoellipticity condition)}.
	\textit{Let us consider the operator \eqref{sumsquares}.
	If {\rm Lie}$\{ X_{0}, X_{1}, \ldots, X_{m} \} (x,t) = \R^{N+1}$ at every $(x,y,t) \in \O$, then $\sum \limits_{i=1}^{m} X_{i}^{2} \, + \, X_{0}$ is hypoelliptic.}

\medskip

We recall that the notation Lie$\{ X_{0}, X_{1}, \ldots, X_{m} \} (x,t)$ denotes the vector space generated 
by the vector fields $\{ X_{0}, X_{1}, \ldots, X_{m} \}$ and their commutator. The \textit{commutator} of two vector fields $W$ and $Z$ acting on $u \in C^{\infty}(\O)$ is defined as $[W,Z] u := WZ u - ZW u$.

As far as we are concerned with the operator $\K_1$ defined in \eqref{K0}, we can write it as follows
\begin{equation}
	\label{protoK}
	\K_1 \, = \,X^{2} \,+ \, Y,
\end{equation}
where
\begin{equation}
	\label{vectorK}
	X = \frac{\p}{\p x} \sim \begin{pmatrix}
						1 \\
						0 \\
						0
			        \end{pmatrix},
	\quad 
	Y = x \frac{\p}{\p y} - \frac{\p}{\p t} \sim \begin{pmatrix}
						0 \\
						x \\
						-1
			        \end{pmatrix}
	\quad \text{and} 
	\quad [X,Y] = \frac{\p}{\p y} \sim \begin{pmatrix}
						0 \\
						1 \\
						0
			        \end{pmatrix}.
\end{equation}
Hence, the vector fields $X, Y$ and $[X,Y]$ form a basis of $\R^{3}$ at every point $(x,y,t) \in \R^{3}$,
so that $K_1$ satisfies the H\"ormander's rank condition.

The commutators are strongly related to a non-Euclidean invariant structure for the Kolmogorov operator, as 
was firstly pointed out by Garofalo and E. Lanconelli in \cite{GL}. Later on, a non commutative algebraic structure was introduced by 
E. Lanconelli and Polidoro in \cite{LP} to replace the Euclidean one in the study of Kolmogorov operators \eqref{K0}. 

\medskip

\noindent 
{\sc Lie group.} \textit{Consider an operator $\K_{1}$ in the form \eqref{K0}. Then $\KK= ( \R^{3}, \bullet)$,}
\begin{equation}
    \label{law-k} 
    (x_{0}, y_{0}, t_{0}) \bullet (x,y,t) = ( x_{0} + x, y_{0} + y - tx_{0} , t_{0} + t).
\end{equation}
\textit{is a group with zero element $(0,0,0)$, and inverse $(x, y, t)^{-1} := (- x, -y - tx, -t)$.}

\medskip 

\noindent
Indeed, if we set $v(x,y,t)=u(x_{0} + x, y_{0} + y - tx_{0} , t_{0} + t)$, then 
\begin{equation*}
	\K_{1} v = 0 \quad \text{if, and only if,} \quad \K_{1} u = 0.
\end{equation*}
Moreover, the operator $\K_{1}$ is invariant with respect to the following family of dilations of $\R^{3}$
\begin{equation}
    \label{dil-k}
    \d_{r}(x,y,t) = ( r x, r^{3} y , r^{2}t) \qquad \text{for every} \; r>0,
\end{equation}
in the sense that if we set $v(x,y,t)=u( r x, r^{3} y , r^{2}t)$, then 
\begin{equation*}
	\K_{1} v = 0 \quad \text{if, and only if,} \quad \K_{1} u = 0.
\end{equation*}

We now introduce a quasi-distance invariant with respect to the group operation ``$\bullet$''. 

\begin{definition}
	\label{dist-def}
	For every $z=(x,y,t), \z=(\x,\y,\t) \in \R^{3}$, we define a quasi-distance $d_{K}(z,\z)$ invariant with 
	respect to the translation group $\KK$ as follows
  	   \begin{equation*}
    		d_{K}(z, \z) = \Big \lvert x - \x \Big \rvert \, + \, \Big \lvert y - \y + (t- \t)\tfrac{x + \x}{2} \Big \rvert^{\tfrac13} \, + \, 
		\Big \lvert t-\t \Big \rvert^{\tfrac12}.
          \end{equation*}
\end{definition}	
Here we recall the meaning of quasi-distance $d_{K}: \R^{3} \times \R^{3} \rightarrow [0, + \infty [$:
	\begin{enumerate}
		\item $d_{K}(z,w) = 0$ if and only if $z=w$ for every $z, w \in \R^{3}$;
		\item $d_{K}(z,w) = d_{K}(w,z)$;
		\item for every $z, w, \z \in \R^{3}$ there exists a constant $C>0$ such that (see Lemma 2.1 of \cite{PDF})
			\begin{equation*}
				d_{K}(z,w) \le d_{K}(z, \z) + d_{K}(\z, w).
			\end{equation*}
	\end{enumerate} 
Moreover, we remark that the quasi-distance $d_{K}$ is homogeneous of degree $1$ with respect to the family of
dilations $\{ \d_{r} \}_{r >0}$ in the sense that for every $z, \z \in \R^{3}$
\begin{equation*}
	d_{K}\left(\d_{r} (z), \d_{r}(\z)\right) \, = \, r \left( d_{K}(z,\z) \right)  \qquad {\rm for \; every} \; r > 0.
\end{equation*}
We are now in position to define the space of H\"older continuous functions $C^{\a}_{K}$.
\begin{definition}
    \label{holdercontinuous}
    Let $\a$ be a positive constant, $\a \le 1$, and let $\O$ be an open subset of $\R^{3}$. We 
    say that a function $f : \O \longrightarrow \R$ is H\"older continuous with exponent $\a$ in $\O$
    with respect to the group $\KK=(\R^{3}, \bullet, \d_{r})$ (in short: H\"older 
    continuous with exponent $\a$, $f \in C^\a_{K} (\O)$) if there exists a positive constant $C>0$ such that 
    \begin{equation*}
        | f(z) - f(\z) | \le C \; d_{K}(z,\zeta)^{\a} \qquad { \rm for \, every \, } z, \z \in \O.
    \end{equation*}
   To every bounded function $f \in C^\a_{K} (\O)$ we associate the norm
    	\begin{equation*}
       		 |f|_{\a, \O} \hspace{1mm} = \hspace{1mm} \sup \limits_\O |f| \hspace{1mm} + \hspace{1mm} 
       		 \sup \limits_{z, \z \in \O \atop  z \ne \z} \frac{|f(z) - f(\z)|}{d_{K}(z,\zeta)^{\a}}.
        \end{equation*}
    Moreover, we say a function $f$ is locally H\"older continuous, and we write $f \in C^{\a}_{K,\loc}(\O)$,
    if $f \in C^{\a}_{K}(\O')$ for every compact subset $\O'$ of $\O$.
\end{definition}
We recall the following result, due to Manfredini (see p. 833 in \cite{M}), where the space $C^{\a}_{K}$ is compared with
the usual Euclidean H\"older space $C^{\a}$.
\begin{proposition}
   \label{manfredini}
   Let $\O$ be a bounded subset of $\R^{3}$.
   If $f \in C^{\a}(\O)$ in the usual Euclidean sense, 
    then $f \in C^{\a}_{K}(\O)$ in the sense of Definition \ref{holdercontinuous}. 
    Vice versa, if $f \in C^\a_{K} (\O)$, then $f \in C^{\b}(\O)$ in the Euclidean sense with $\b = \tfrac{\a}{3}$.
\end{proposition}

We remark that local H\"older regularity assumption we assume on the coefficients of the operator $\K$
is less restrictive than the global H\"older regularity, as pointed out by Pascucci and Pesce (see Example 1.3, \cite{PaPe}). 
Indeed, for every $y, \y, t, \t \in \R$ with $t \ne \t$,
let us consider the following couple of points in $\R^{3}$
 \begin{equation}
 	\label{punti}
 	z = \left( \frac{\y - y}{t - \t}, y,  t \right) \qquad \text{and} \qquad \z = \left(\frac{\y - y}{t - \t}, \y, \t \right),
 \end{equation}
then we have $d ( z, \z) \, = \, \abs{t-\t}^{\tfrac12}$.  
Since $y$ and $\y$ are arbitrary real numbers, we see that points in $\R^{3}$ that are far from each other in the Euclidean sense 
can be very close with respect to the distance $d$. It follows that, if a function $f(x,y,t) = f(y)$ depends 
only on $y$ and it belongs to $C^{\a}_{K} (\R^{3})$ for some positive $\a$, then it must be constant. In fact, for $z, \z$ 
as defined in \eqref{punti}, we have
 \begin{equation*}
 	\abs{f(y) - f(\y)} \, = \, \abs{f(z) - f (\z) } \, \le \, C \abs{t-\t}^{\a}
 \end{equation*}
 for some positive constants $C$, $\a$ and for any $y,\y \in \R$ and $t \ne \t$.

\subsection{H\"older continuous coefficients}
In this work we consider classical solutions to the equation $\K u = f$. 
We first recall the notion of Lie derivative $Y u$ of a function $u$ with respect to the vector field $Y$ defined in \eqref{vectorK}: 
\begin{equation} \label{lie-diff}
		Yu(x,y,t) := \lim \limits_{s \rightarrow 0} \frac{u(\g(s)) - u(\g(0))}{s}, \qquad \g(s) = (x, y + s x, t - s).
\end{equation}
Note that $\g$ is the integral curve of $Y$, i.e. $\dot \g (s) = Y(\g(s))$. 
Clearly, if $u \in C^{1}(\O)$, with $\O$ open subset of $\R^{3}$, then $Y u (x,y,t)$ agrees with $x \p_{y} u(x,y,t) - \p_{t} u (x,y,t)$ 
considered as a linear combination of the derivatives of $u$. 

\begin{definition}
	\label{solution}
	A function $u$ is a solution to the equation $\K u = f$ in a domain $\O$ of $\R^{3}$
	if the derivatives $\p_{x} u, \p_{x}^{2} u$ and the Lie derivative $Yu$ exist as continuous functions in $\O$, and the 
	equation
	\begin{equation*}
	     \frac{\p}{\p x} \biggl(a (x,y,t) \frac{\p u}{\p x} \biggr) \, + \, b(x,y,t) \, \frac{\p u}{\p x} \, + \,  x \frac{\p u}{\p y} \, 
	     - r(x,y,t) u - \, \frac{\p u}{\p t} = f(x,y,t)
	\end{equation*} 
	is satisfied at any point $(x,y,t) \in \O$. Moreover, 
	we say that $u$ is a \emph{classical supersolution}  to $\K u = 0$ if $f \le 0$ in $\O$, and we write $\K u \le 0$.
	We say that $u$ is a \emph{classical subsolution} if $-u$ is a classical supersolution.
\end{definition}
A fundamental tool in the classical regularity theory for Partial Differential Equations are the Schauder estimates.
In particular, we recall the result proved by Manfredini in \cite{M} (see Theorem 1.4)  for classical solutions to 
$\K u = f$, where the natural functional setting is
 \begin{equation*}
	C^{2 + \a}_{K} (\O) = \left\{ u \in C^{\a}_{K} (\O) \; \mid \; \p_{x}u , \p^2_{x} u, Y u \in C^{\a}_{K} (\O) \right\},
   \end{equation*}
and $C^{\a}_{K} (\O)$ is given in Definition \ref{holdercontinuous}. Moreover, if $u \in C^{2 + \a}_{K} (\O)$ then we define the norm
\begin{equation*}
	| u |_{2 + \a, \O } := | u |_{\a, \O } \; + \; | \p_{x}u |_{\a, \O } \; + \;  |\p^2_{x}u |_{\a, \O }  \; + \;  |Y u |_{\a, \O }.
   \end{equation*}
Clearly, the definition of $C^{2+\a}_{K,\loc} (\O)$ follows straightforwardly from the definition of $C^{\a}_{K,\loc} (\O)$. 
We refer to the survey paper \cite{APsurvey} for a more recent bibliography on this subject.
\begin{theorem}
    \label{schauder1}
    Let us consider an operator $\K$ of the type \eqref{K} satisfying assumptions \textbf{(H$_{\K}$)} with $0 < \a \le 1$. 
    Let $\O$ be an open subset of $\R^{3}$, $f \in C^{\a}_{K,\loc}(\O)$ and let $u$ be a classical solution to 
    $\K u = f$ in $\O$. Then for every $\O^{'} \subset \subset \O^{''} \subset \subset \O$ 
    there exists a positive constant $C$ such that
    \begin{equation*}
        \label{scauder2}
        | u |_{2 + \a, \O^{'} } \le C \Big( \sup \nolimits_{\O^{''}} |u| \; + \; |f|_{\a, \O^{''} } \Big).
    \end{equation*}
\end{theorem}

We also recall the notion of fundamental solution. 
\begin{definition}
	\label{funsolk}
	We say a function $\GK: \R^{3} \times \R^{3} \rightarrow \R$ is a fundamental solution for $\K$ if
	 \begin{enumerate}
	 	\item for every $(x_{0}, y_{0}, t_{0}) \in \R^{3}$ the function $x \mapsto \GK (x, y, t; x_{0}, y_{0}, t_{0})$:
			\begin{enumerate} 
				\item belongs to $L^{1}_{\loc} (\R^{3})$; \\
				\item is a classical solution of $\K u = 0$ in $\R^{3} \setminus \{ (x_{0}, y_{0}, t_{0}) ) \}$
					in the sense of Definition \ref{solution};
			\end{enumerate}
		\item for every bounded function $\phi \in C(\R^{2})$, we have that
			$$
 				 u(x,y,t)=\int_{\R^{2}} \GK(x,y,t;\x,\y,t_0) \, \phi(\x, \y) \, d\x \, d\y, 
			$$
			is a classical solution of the Cauchy problem \eqref{cauchyproblem};
		\item The function $\GK^{\! \! \! \!*}(x, y, t; x_{0}, y_{0}, t_{0}) := \GK (x_{0}, y_{0}, t_{0}; x, y, t)$ satisfies 1. and 2. with $\K$
			replaced by its adjoint operator $\K^{*}$ as defined in \eqref{Kadj}.
	\end{enumerate}
\end{definition}

The existence of a fundamental solution $\GK$ for the operator $\K$ has widely been investigated over the years, and as we have already pointed out in the Introduction of this paper the Levy parametrix method provides us with a fundamental solution for the operator $\K$ under global H\"older regularity assumptions for the coefficients. Among the first results of this type we recall \cite{Weber}, \cite{IlIn} and \cite{Sonin}. We summarize here the main results of the articles \cite{P2}, \cite{PDF} and \cite{LPP} on the existence and bounds for the fundamental solution under the following assumption for the coefficients of the operator $\K$:
\begin{itemize} 
	\item[] $a,\frac{\p a}{\p x},b,r \in C^{\a}_{K}(\R^{3})$ and there exist two positive constants $\l, \Lambda$ such that 
	\begin{align}
		\label{Hmix}
		\l \le a(x,y,t) \le \Lambda \quad \left| \tfrac{\p a}{\p x} \right|, \abs{b(x,y,t)}, \abs{r(x,y,t)} \le \Lambda  \quad
		\text{for every} \quad (x,y,t) \in \R^3.
	\end{align}
\end{itemize}
For more references on this subject we refer to the survey paper \cite{APsurvey}.

\begin{theorem}
	\label{ex-prop}
	Let $\K$ be an operator of the form \eqref{K} under the assumption \eqref{Hmix}. Then there exists a 
	fundamental solution $\GK$ in the sense of Definition \ref{funsolk}. Moreover, for every $(x_{0}, y_{0}, t_{0}) \in \R^{3}$,
	$\GK$ belongs to $C^{2 + \a}_{\loc}(\R^{3} \setminus \{ (x_{0}, y_{0}, t_{0}) ) \}$ and
	the following properties hold:
	\begin{enumerate}
	\item Support of $\GK$: for every $(x,y,t), (\x, \y, \t) \in \R^{3}$ with $t \le \t$
		\begin{equation*}
			\GK (x,y,t; \x, \y, \t) = 0;
		\end{equation*}
	\item Reproduction property: for every $(x,y,t), (x_{0}, y_{0}, t_{0}) \in \R^{3}$ and $\tau \in \R$ with $t_{0} < \t < t$: 
		\begin{equation*}
			\GK(x,y,t; x_{0}, y_{0},t_{0}) = \int \limits_{\R^{2}} \GK(x,y,t;\x, \y, \t) \, \GK(\x, \y, \t; x_{0}, y_{0}, t_{0}) \, d\x \, d\y;
		\end{equation*}
	\item Integral of $\GK$: the following bound holds true
		\begin{equation}
			\label{integralK}
             e^{-\Lambda (t - \tau)} \, \le \, 
 			\int \limits_{\R^{2}} \GK(x,y,t; \x, \y, \t) \; d \x \, d\y 
 			\, \le \, e^{\Lambda (t - \tau)};
 		\end{equation} 
	\item let $I=]T_{0}, T_{1}[$ be a bounded interval, then there exist four positive constants $\l^{+}$, $\l^{-}$,
			$C^{+}$, $C^{-}$ such that for every $(x,y,t), (\x, \y, \t) \in \R^{3}$ with $T_{0} < t < \t < T_{1}$
			\begin{align}
				\label{boundK}
				C^{-} \, \GK^{\! \! \! \! \l^{-}} (x,y,t; \x, \y, \t) \, \le \, \GK(x,y,t; \x, \y, \t) \, \le \, 
				C^{+} \, \GK^{\! \! \! \! \l^{+}} (x,y,t; \x, \y, \t).
			\end{align}
		The constants $\l^{+}$, $\l^{-}$, $C^{+}$, $C^{-}$ depend only on $\K$ and $T_{1}-T_{0}$. 
		$\GK^{\! \! \! \! \l^{-}}$ and $\GK^{\! \! \! \! \l^{+}}$ respectively denote the fundamental solution of 
		$\K_{\l^{-}}$ and $\K_{\l^{+}}$, defined in \eqref{Glambda} and \eqref{Klambda} respectively.	
	\end{enumerate}
	Furthermore, for every $(x_{0}, y_{0}, t_{0}) \in \R^{3}$ also
	$\GK^{*}$ belongs to $C^{2 + \a}_{\loc}(\R^{3} \setminus \{ (x_{0}, y_{0}, t_{0}) ) \}$.
\end{theorem}

The properties \emph{1.} and \emph{2.} of the above statement have been proved in \cite{P2} and \cite{PDF}. The inequalities \eqref{integralK} follow from the comapison principle for classical solutions, as the functions $e^{-\Lambda (t - \tau)}$ and $e^{\Lambda (t - \tau)}$ are respectively subsolution and supersolution to the Cauchy problem \eqref{cauchyproblem} with inital datum $\phi(x,y)= 1$. We remark that the constants $\l^{+}$, $\l^{-}$, $C^{+}$, $C^{-}$ appearing in the bounds \eqref{boundK} proved in \cite{P2, PDF} also depend on the $C^\alpha(\R^N \times I)$ norm of the coefficients $a,\frac{\p a}{\p x},b,r$. We rely here on the bounds proved in \cite{LPP}, where the dependence on the regularity of the coefficients is removed thanks to the Harnack inequality proved by Golse, Imbert, Mouhot and Vasseur in \cite{GIMV}. We conclude this section with the following Gaussian bound for $\GK$. 
\begin{corollary} \label{corollary}
	Let $(x,y,t) \in \R^{3}$, with $t > t_{0}$. Then there exist two positive constants $\overline C$, only depending on the operator $\K$, and $R_{0}$, also depending on $(x,y,t)$, on $t_0$, such that
	\begin{equation*}
		\GK \left( x, y, t; \x, \y, \t \right) \le \overline C e^{-\overline C \, \frac{\x^{2} + \y^{2}}{t - \t}},
	\end{equation*}
for every $(\x, \y) \in \R^{2}$ such that $\left\| (\x, \y) \right\| \ge R_{0}$ and for every $\t \in \R$ with $t_{0} < \t < t$. 
\end{corollary}
The proof of this result directly follows from the upper bound \eqref{boundK} combined with the explicit expression of the fundamental solution \eqref{Glambda}. We refer to the Lemma 3.1 of \cite{Punic} for the proof, that will be omitted here.

\setcounter{equation}{0} 
\section{Geometric setting and fundamental solution for $\L$}
\label{preliminaries}
The aim of this section is to recall some known properties and results concerning the operator $\L$. Moreover, we define the space of H\"older continuous functions $C^{\a}_{L}$ associated to $\L$ and we compare it with the space 
of H\"older continuous functions $C^{\a}_{K}$ defined in Definition \ref{holdercontinuous}.

\subsection{The operator $\L_{1}$}
Let us consider the operator $\L_{1}$ introduced in \eqref{Llambda} as the prototype operator for the operator $\L$. 
As we have already pointed out in the Introduction of this paper, the function $\GL^{\! \! 1}$ defined in \eqref{fun-yor} 
is the fundamental solution $\GL^{\! \! 1}$ of $\L_{1}$. Its expression agrees with that of the density of the process $(W_{t}, A_{t})_{t \ge 0}$ in \eqref{processY}, first considered by Yor in \cite{Y}.

As far as we are concerned with the invariance properties of $\L_{1}$, Monti and Pascucci observe in \cite{MP} that it 
is invariant with respect to the group operation
on $\R^{+} \times \R^{2}$:
\begin{equation} \label{trasl}
	(\x, \y, \t) \circ (x,y,t) = (\x \, x , \y + \x \,y , \t + t) .
\end{equation}
Indeed, if we set $v(x,y,t)=u(\x \, x , \y + \x \,y , \t + t)$, then $\L_{1} v = 0$ if, and only if, $\L_{1} u = 0$.
We also remark that 
\begin{equation*}
	\GG := (\R^{+} \times \R^{2} , \circ) 
\end{equation*}
is a Lie group, where the identity $\I_{\GG}$ and the inverse of the element $(x,y,t)$ are defined as 
\begin{equation*}
	\I_{\GG} = (1,0,0), \qquad (x,y,t)^{-1} = \left(\tfrac{1}{x}, -\tfrac{y}{x}, - t\right).
\end{equation*}
Let us notice that the translation defined in \eqref{trasl} reflects the mixed nature of the stochastic process $(S_t, A_t)_{t \ge 0}$ defined in \eqref{processY}. Indeed its first component $S_t$ is log-normal, then is related to a multiplicative group, while the component $A_t$ is defined as the integral of $S_t$, then is related to an additive group. In particular, the null element of the group is $(1,0,0)$, the left-translation $(r,0,0) \circ (x,y,t)$ acts as a dilation with respect to $(x,y)$, while the left-translation $(1,\eta,t) \circ (x,y,t)$ acts as an Euclidean translation with respect to $(y,t)$
\begin{equation*}
	(r,0,0) \circ (x,y,t) = (rx , ry , t), \qquad (1,\eta,t) \circ (x,y,t) = (x, \y + \,y , \t + t).
\end{equation*}

As far as we are concerned with the regularizing properties of the operator $\L_{1}$, one can easily prove it is 
hypoelliptic in the sense of \eqref{hypo}. Indeed, we can write the vector fields associated to $\L_{1}$ as follows
\begin{align}
	\label{vec0}
	X = x \p_{x} \sim & \begin{pmatrix}
           						x \\
           						0 \\
           						0
        						 \end{pmatrix}, \quad
	Y = x \p_{y} - \p_{t} \sim \begin{pmatrix}
           							0 \\
           							x \\
           							-1
        						 	   \end{pmatrix}, \quad \text{and} \quad
	[X, Y] = x \p_{y} \sim \begin{pmatrix}
           							0 \\
           							x \\
           							0
        						 	   \end{pmatrix}.			   									
\end{align}
Thus, Lie$\{ X, Y, [X,Y] \} (x,y,t) = \R^{+} \times \R^{2}$ for every $(x,y,t) \in \R^{+} \times \R^{2}$, hence $\L_{1}$
satisfies the H\"ormander's hypoellipticity condition.

\subsection{The optimal control problem for $\L_{1}$}
We now introduce the function $\psi$ appearing in the formula \eqref{boundL}. Let us consider the vector fields $X$ and $Y$ defined in \eqref{vec0} associated to the operator $\L_{1}$. We consider the following optimal control problem. For any \emph{end point} 
$(x_{0}, y_{0}, t_{0}) \in \R^{+} \times \R^{2}$ and \emph{starting point} $(x_{1}, y_{1}, t_{1}) \in \R^{+} \times \R^{2}$,
with $t_{1} > t_{0}$:
\begin{align}
	\label{optcont}
	\psi(x_{1},y_{1},t_{1}; x_{0},y_{0},t_{0}) \, := \, &\inf \limits_{\w \in L^{1}([0,T])} \, 
	\int \limits_{0}^{T} \w^{2} (\t) \, d\t 
	\qquad \text{\rm subject to constraint} \\ \nonumber
	&\begin{cases}
		\dot{x}(s) = \w(s) \, x(s) \\
		\dot{y}(s) = x(s)  \qquad \qquad  0 \le s \le T\\
		\dot{t}(s)  = -1
	\end{cases} \\ \nonumber
	(x,y,t)(0) = & \,(x_{1}, y_{1}, t_{1}), \qquad (x,y,t)(T) = (x_{0}, y_{0}, t_{0}).
\end{align}
The constraint $\dot{t}(s)  = -1$ implies that admissible paths satisfy $t(s) = t_{1} - s$, hence $T=t_{1} - t_{0}$.
For this reason, in the sequel we drop the time variable, and we set $T := t_{1} - t_{0}$. Moreover, as $x(s) >0$ for every $s$, the second equation yields $y_1 < y_0$. 
The knowledge of the explicit expression of the function $\psi$ is particularly important, and we summarize here some quantitative informations about it in terms of the following function
\begin{equation*}
	g(r) = \begin{cases}
			\frac{\sinh (\sqrt{r})}{\sqrt(r)}, \qquad &r > 0, \\
			1 \quad &r=0,\\
			 \frac{\sinh (\sqrt{-r})}{\sqrt(-r)}, \qquad &- \pi^{2} < r < 0.
		 \end{cases}	
\end{equation*}
Indeed, for every $(x_{1}, y_{1}, t_{1}), (x_{0}, y_{0}, t_{0}) \in \R^{+} \times \R^{2}$, with $t_{0} < t_{1}$
and $y_{0} > y$, we have 
\begin{equation}
	\label{psi-open}
	\begin{cases}
		\psi (x_{1}, y_{1}, t_{1} ; x_{0}, y_{0}, t_{0}) = E(t_{1} - t_{0}) + \frac{4(x_{1} + x_{0})}{y_{0} - y_{1}} 
		- 4 \sqrt{E + \frac{4x_{1}x_{0}}{(y_{0} - y_{1})^{2}}}, \quad 
		\text{if} \, \, E \ge - \frac{\pi^{2}}{(t_{1} - t_{0})^{2}}& \! \! \! \!; \\
		\\
		 \psi (x_{1}, y_{1}, t_{1} ; x_{0}, y_{0}, t_{0}) = E(t_{1} - t_{0}) + \frac{4(x_{1} + x_{0})}{y_{0} - y_{1}} 
		+ 4 \sqrt{E + \frac{4x_{1}x_{0}}{(y_{0} - y_{1})^{2}}},  \\
		\hfill \text{if} \, \, - \frac{4\pi^{2}}{(t_{1} - t_{0})^{2}} < E < - \frac{\pi^{2}}{(t_{1} - t_{0})^{2}} & \! \! \! \!,
	\end{cases}
\end{equation}
where 
\begin{equation*}
	E = \frac{4}{(t_{1} - t_{0})^{2}} \, g^{-1} \left( \frac{y_{0} - y_{1}}{(t_{1} - t_{0}) \sqrt{x_{1} x_{0}}} \right).
\end{equation*}
For further informations we refer to \cite{CPR}, Section 4, where also the solution of the control problem \eqref{optcont} 
is provided. Moreover, we recall that one of the results of \cite{CPR} are the following bounds for the fundamental 
solution $\GL^{\! \! 1}$ constructed by Geman and Yor in \cite{GemanYor}:
\begin{align}
	\label{bound}
	\frac{c_{\e}^{-}}{x_{0}^{2} (t - t_{0})^{2}} \exp & \left( -C^{-} \psi (x,y + x_{0} \e (t-t_{0}) , t - \e (t-t_{0}); 
	x_{0}, y_{0},t_{0}) \right) \le \\ \nonumber
	& \qquad \le \GL^{\!\! 1} (x,y,t; x_{0}, y_{0},t_{0}) \le \frac{C_{\e}^{+}}{x_{0}^{2} (t - t_{0})^{2}} \exp \left(-c^{+} \psi (x,y - x_{0} \e, t + \e; 
	x_{0}, y_{0},t_{0}) \right),
\end{align}
where $\psi$ is the cost function of the optimal control problem \eqref{psi}.

As the vector fields $X = x \partial_x$ and $Y=x \partial_y - \partial_t$ are invariant with respect to the left translation ``$\circ$'' in \eqref{trasl}, it turns out that also the solution $\psi$ to the optimal control problem \eqref{optcont} is invariant with respect to  $\GG = (\R^{+} \times \R^{2}, \circ)$. In particular we have
\begin{equation*}
	\psi (x_{1}, y_{1}, t_{1}; x_{0}, y_{0}, t_{0}) = \psi ( (x_{0}, y_{0}, t_{0})^{-1} \circ (x_{1}, y_{1}, t_{1}); 1,0,0).
\end{equation*}
Hence, from now on we fix the final condition $(x_{0}, y_{0}, t_{0})=(1,0,0)$ in the optimal control problem \eqref{optcont}, 
and then use the invariance property to solve it with a general initial condition, and we use the simplified notation $\psi(x,y,t) = \psi(x,y,t; 1,0,0)$.
For all of the above reasons, the optimal control problem \eqref{optcont} now reads as follows for a general starting point 
$(x,y,t) \in \R^{+} \times \R^{-} \times \R^{+}$:
\begin{align}
	\label{psi}
	\psi(x,y,t) \, &:= \, \inf \limits_{\w \in L^{1}([0,t])} \, \int \limits_{0}^{t} \w^{2} (\t) \, d\t 
	\qquad \text{\rm subject to constraint} \\ \nonumber
	&
	\begin{cases}
		\dot{x}(s) = \w(s) \, x(s), \quad x(0) \! \! \! \!&= x, \quad x(t) = 1, \\
		\dot{y}(s) = x(s), \quad  \hfill    y(0) \! \! \! \!&= y, \; \quad y(t) = 0.
	\end{cases} 
\end{align}

\subsection{The space $C^{\a}_{L}$}
In order to define the H\"older spaces relevant to the operator $\L$ we note that the operators $\L$ and $\K$ agree in every compact set of $\R^+ \times \R^2$. We then borrow the regularity theory developed for the opeator $\K$, and described in the above subsection, to obtain analogous results for the operator $\L$. This point of view was adopted in the work \cite{CPR} to obtain an invariant Harnack inequality for $\L$. 

Consider a function $f=f(x,y,t)$ defined in $\R^{+} \times \R^{2}$, and let $(\xi, \eta, \tau)$ be a point
in $\R^{+} \times \R^{2}$. In accordance with the operation \eqref{trasl}, we define the function 
\begin{equation*}
 \widetilde f(x,y,t) := f (\x \, x , \y + \x \,y , \t + t).
\end{equation*}
We note that
\begin{equation*}
 f(x,y,t) = \widetilde  f \left( \tfrac{x}{\xi}, \tfrac{y - \eta}{\xi}, t - \tau \right),
\end{equation*}
and we apply the Definition \ref{holdercontinuous} to $\widetilde f(x,y,t)$ in a neighborhood of $(1,0,0)$. We find 
\begin{equation}
	\label{eq1}
\begin{split}
 |f(x,y,t) - f (\x,\y,\t)| = & \left| \widetilde f\left( \tfrac{x}{\xi}, 
 \tfrac{y - \eta}{\xi}, t - \tau \right) - \widetilde f(1,0,0) \right| \le \\
 & C \left( \left| \tfrac{x - \x}{\xi} \right| + \left| \tfrac{y - \eta}{\xi} +
 (t-\tau) \tfrac{x + \x}{2\xi} \right|^{1/3} + |t-\tau|^{1/2} \right)^\alpha.
\end{split}
\end{equation}
Let us remark that the operators $\L$ and $\K$ are comparable only when the points $x$ and $\xi$ are close each other. Indeed,  if we exchange the role of the points $(x,y,t)$ and $(\x, \y, \t)$, we find the inequality 
\begin{equation*}
 |f(\x,\y,\t) - f (x,y,t)| \le C \left( \left| \tfrac{\x - x}{x} \right| + \left| \tfrac{\y - y}{x} +
 (\t-t) \tfrac{\x + x}{2x} \right|^{1/3} + |t-\tau|^{1/2} \right)^\alpha,
\end{equation*} 
which doesn't agree with \eqref{eq1}, unless $\xi$ and $x$ have similar size. For this reason, we give the following definition of quasi-distance and H\"older continuous function with respect to the operation ``$\circ$'', which is equivalent to \eqref{eq1} when $\frac{x}{\xi}$ is close to $1$.


\begin{definition}
	\label{dL}
	For every $z=(x,y,t), \z=(\x,\y,\t) \in \R^{+} \times \R^{2}$, we define a symmetric quasi-distance $d_{L}(z,\z)$ 
	invariant with respect to the translation group $\GG$ as follows
  	   \begin{equation*}
    		d_{L}(z, \z) = \left| \tfrac{x - \x}{\sqrt{x \xi}} \right| + 
    		\left| \tfrac{y - \eta + (t-\tau) \tfrac{x + \x}{2}}{\sqrt{x \xi}} \right|^{1/3} + |t-\tau|^{1/2}.
          \end{equation*}

\end{definition}
\begin{definition}
    \label{holdercontinuousL}
    Let $\a$ be a positive constant, $\a \le 1$, and let $\O$ be an open subset of $\R^{+} \times \R^{2}$. We 
    say a function $f : \O \longrightarrow \R$ is H\"older continuous with exponent $\a$ in $\O$
    with respect to the group $\GG=(\R^{+} \times \R^{2}, \circ)$ (in short: H\"older 
    continuous with exponent $\a$, $f \in C^\a_{L} (\O)$) if there exists a positive constant $C>0$ such that 
    \begin{equation}\label{e-assumption2}
 		 | f(x,y,t) - f (\xi, \eta, \tau) |  \le  C d_{L}(z, \z)^\alpha,
	\end{equation}
	for every $(x,y,t), (\xi, \eta, \tau) \in \O$. 
    Moreover, we say a function $f$ is locally H\"older continuous, and we write $f \in C^{\a}_{L,\loc}(\O)$,
    if $f \in C^{\a}_{L}(\O')$ for every compact subset $\O'$ of $\O$.
\end{definition}

As the definitions $C^{\a}_{L}(\O')$ and $C^{\a}_{K}(\O')$ agree in every compact subset $\O'$ of $\R^{+} \times \R^{2}$, the following statement is an immediate consequence of Proposition \ref{manfredini}. For this reason, we omit its proof, which is immediate.

\begin{proposition}
   \label{coincide}
   Let $\O'$ be a compact subset of $\R^{+} \times \R^{2}$. If $f \in C^{\a}(\O)$ in the usual Euclidean sense, 
    then $f \in C^{\a}_{L}(\O')$ in the sense of Definition \ref{holdercontinuousL}. 
    Vice versa, if $f \in C^\a_{L} (\O')$, then $f \in C^{\b}(\O')$ in the Euclidean sense with $\b = \tfrac{\a}{3}$.
\end{proposition}

\begin{definition}
	\label{solutionL}
	A function $u$ is a solution to the equation $\L u = f$ in a domain $\O$ of $\R^{+} \times \R^{2}$
	if the derivatives $x \p_{x} u, x^{2} \p_{x}^{2} u, $ and the Lie derivative $Yu$ exist as continuous functions in 
	$\O$, and the equation
	\begin{equation*}
	     x \, \frac{\p u}{\p x} \left( a (x,y,t) \, x \, \frac{\p u}{\p x} \right) \, + \, b(x,y,t) \, x \, \frac{\p u}{\p x}  \, + \, x \, 
	     \frac{\p u}{\p y} - r(x,y,t) u - \, \frac{\p u}{\p t}, = f(x,y,t)
	\end{equation*} 
	is satisfied at any point $(x,y,t) \in \O$.
\end{definition}

 \begin{equation*}
	C^{2 + \a}_{L} (\O) = \left\{ u \in C^{\a}_{L} (\O) \; \mid \; x \p_{x}u , x^{2} \p^2_{x} u, Y u \in C^{\a}_{L} (\O) \right\},
   \end{equation*}
where $C^{\a}_{L} (\O)$ is given in Definition \ref{holdercontinuousL}. 
Clearly, the definition of $C^{2+\a}_{L,\loc} (\O)$ follows straightforwardly from the definition of $C^{\a}_{L,\loc} (\O)$. 

\subsection{The fundamental solution $\GL$}
We now focus on the fundamental solution $\GL$ for the operator $\L$. 
\begin{definition}
	\label{funsoll}
	We say a function $\GL: \left( \R^{+} \times \R^{2} \right) \times \left( \R^{+} \times \R^{2} \right) \rightarrow \R$ is a fundamental solution 
	for $\L$ if
	 \begin{enumerate}
	 	\item for every $(x_{0}, y_{0}, t_{0}) \in \R^{+} \times \R^{2}$ the function 
			$x \mapsto \GL (x, y, t; x_{0}, y_{0}, t_{0})$:
			\begin{enumerate} 
				\item belongs to $L^{1}_{\loc} (\R^{+} \times \R^{2}) \cap C^{\infty}(\R^{+} \times \R^{2} 
					\setminus \{ (x_{0}, y_{0}, t_{0}) ) \}$; \\
				\item is a classical solution of $\L u = 0$ in $\R^{+} \times \R^{2} \setminus 
					\{ (x_{0}, y_{0}, t_{0}) ) \}$ in the sense of Definition \ref{solutionL};
			\end{enumerate}
		\item for every bounded function $\phi \in C(\R^{2})$, we have that
		\begin{equation} \label{rep2}
		 u(x,y,t)=\int_{\R^{2}} \GL(x,y,t;\x,\y,0) \, \phi(\x, \y) \, d\x \, d\y, 
		\end{equation}
			is a classical solution of the Cauchy problem
			\begin{equation} \label{cauchyproblemL}
			\begin{cases}
    				\L u = 0,  \qquad &(x,y,t) \in \R^{+} \times \R \times \R^{+} \\
    				u(x,y,0)=\varphi (x,y) \quad&(x,y) \in \mathbb{R}^2.
 			\end{cases}
		\end{equation}
		\item The function $\GL^{\! \! *}(x, y, t; x_{0}, y_{0}, t_{0}) := \GL (x_{0}, y_{0}, t_{0}; x, y, t)$ satisfies 1. 
			and 2. with $\L$ replaced by its adjoint operator $\L^{*}$ as defined in \eqref{L*}.
	\end{enumerate}
\end{definition}
Under the following assumption \eqref{HmixL}, which is stronger than \textbf{(H$_L$)}, Cibelli, Polidoro and Rossi prove the existence of the fundamental solution $\GL$ for $\L$ and bounds analogous to \eqref{bound} by applying methods coming from the stochastic theory
(see Proposition 3.7 and Theorem 1.3 of \cite{CPR}, respectively). We summarize here the main results of the paper \cite{CPR}, under the following assumption for the coefficients $a$ and $b$:

\begin{itemize}
 	\item[] $a$,$b$ $\in C^{\infty}(\R^{+} \times \R^2)$. Moreover, there exist two positive constants $\l$, $\Lambda$ 		such that 
		\begin{align}
			\label{HmixL}
			& \lambda \le a(x,y,t) \le \Lambda, \quad \left| \tfrac{\p}{\p x}  \left( x \, a(x,y,t) \right)\right|,
			\left| \tfrac{\p}{\p x}  \left( x \, a(x,y,t) \right)\right|\le \Lambda, 
		\end{align} 
for every $(x,y,t) \in \R^{+} \times \R^{2}$.
\end{itemize}

\begin{theorem}
	\label{ex-propL}
	 Let $\L$ be an operator of the form \eqref{Ldiv} under the assumption \eqref{HmixL}.  
	 Then there exists a fundamental solution $\GL$ in the sense of Definition \ref{funsoll}. Moreover, the following 
	 properties hold:
	\begin{enumerate}
		\item Support of $\GL$: for every $(x,y,t)$, $(\x, \eta, \t) \in \R^{+} \times \R \times ]0,T]$ with $t \le \t$
			and $y \ge \y$
			\begin{equation*}
				\GL (x,y,t; \x, \y, \t) = 0;
			\end{equation*}
		\item Reproduction property:  for every $(x,y,t)$, $(x_{0}, y_{0}, t_{0})$, $(\x, \eta, \t) \in \R^{+} \times \R 
		\times ]0,T]$ 
			with $T \le t > \t > t_{0} > 0$ 
			\begin{equation*}
				\GL (x,y,t; x_{0}, y_{0}, t_{0}) = \int \limits_{\R^{+} \times \R}
				\GL (x,y,t; \x, \y, \t) \, \GL(\x, \y, \t; x_{0}, y_{0}, t_{0}) \, d\x \, d\y;
			\end{equation*}
		\item Integral of $\GL$: there exists a positive constant $\overline{C}$ depending on $t-\t$ and such that 
			$\overline{C} \rightarrow 1$ as $t \rightarrow \t$ for which
			\begin{equation*}
				\int \limits_{\R^{+} \times \R} \GL (x,y,t; \x, \y, \t) \, d\x \, d\y  \, = \, 1,     \qquad 
				\int \limits_{\R^{+} \times \R} \GL (x,y,t; \x, \y, \t) \, dx \, dy  \, = \, \overline{C};
			\end{equation*}
		\item Bounds for $\GL$: for every $\e \in ]0,1[$, and $T>0$ there exist two positive constants $c_{\e}^{-}, 
			C_{\e}^{+}$ depending on $\e$, on $T$ and on the operator $\L$, and two positive constants $C^{-}, 
			c^{+}$, only depending on the operator $\L$, such that 
			\begin{align}
				\label{boundL}
				\frac{c_{\e}^{-}}{x_{0}^{2} (t - t_{0})^{2}} & \exp \left( -C^{-} \psi (x,y + x_{0} \e (t-t_{0}) ,
				 t - \e (t-t_{0}); x_{0}, y_{0},t_{0}) \right) \le \\ \nonumber
				& \qquad \le \GL (x,y,t; x_{0}, y_{0},t_{0}) \le \frac{C_{\e}^{+}}{x_{0}^{2} (t - t_{0})^{2}} \exp 
				\left(-c^{+} \psi (x,y - x_{0} \e, t + \e; x_{0}, y_{0},t_{0}) \right)
			\end{align}
			for every $(x,y,t) \in \R^{+} \times ]-\infty, y_{0} - x_{0}\e(t-t_{0})[ \times ]t_{0}, T]$, where $\psi$ is the 
			value function for the optimal control problem \eqref{psi}.
	\end{enumerate}
\end{theorem}
We remark that the bounds obtained in \eqref{boundL} for the operator $\L$ by \cite{CPR} are analogous to the bounds \eqref{boundK} 
	obtained for the fundamental solution of the Kolmogorov operator $\K$. Let us consider the fundamental solutions $\GL^{\pm}$ of the operators 
	\begin{equation}
		\label{piumenoL}
		\L^{\pm} u = \l^{\pm} x^{2} \frac{\p^{2}u}{\p x^{2}} \, + \, x \frac{\p u}{\p x} \, + \, x \frac{\p u}{\p y}  \, - \, \frac{\p u}{\p t} ,
		\quad (x,y,t) \in R^{+} \times \R \times \R^{+}.
	\end{equation}
By applying the bounds \eqref{boundL} to $\GL$ and to $\GL^{\pm}$, we obtain the following corollary to the previous Theorem \ref{ex-propL} 
	 (see Proposition 1.5 of \cite{CPR}). 
	\begin{corollary}
		\label{corollaryL}
		For every $\e \in ]0,1[$, there exist the fundamental solutions $\G^{\pm}$ of the operators \eqref{piumenoL}, 
		and positive constants $k^{\pm}$ such that
		\begin{align*}
			k^{-} \GL^{-} (x, y + x_{0} \e (t &- t_{0} + 1), t- \e (t - t_{0} + 1); x_{0}, y_{0}, t_{0}) \le \\
				\le \, \GL (x,&y,t ; x_{0}, y_{0}, t_{0}) \le \\
				&\le \, k^{+} \GL^{+} \left(x, y - x_{0} \frac{\e}{1 - \e} (t - t_{0} + 1), t + \frac{\e}{1 - \e} (t - t_{0} + 1); x_{0}, y_{0}, t_{0}\right)
		\end{align*}
		for every $(x,y,t), (x_{0}, y_{0}, t_{0}) \in \R^{+} \times \R \times ]0, T]$, with $y + x_{0} \e (t - t_{0} + 1) < y_{0}$ and 
		$t > t_{0} + \e / (1 - \e)$.
	\end{corollary}

\setcounter{equation}{0} 
\section{Main results}
\label{proof}
This section is devoted to the proof of our main results on the existence of the fundamental solution for the operators 
$\K$ and $\L$. Our approach relies on the local regularity properties of the solutions and on the bounds for the fundamental solution. 

Let us consider first the operator $\K$. We build a sequence of operators $\left( \K_n \right)_{n \in \N}$ satisfying the 
hypotheses of Theorem \ref{ex-prop}, then a fundamental solution $\GK^{\! \! \! n}$ exists for every $n \in \N$. Moreover, 
the sequence $\left( \GK^{\! \! \! n} \right)_{n \in \N}$ is equibounded by \eqref{boundK}, and locally equicontinuous, thanks 
to the Schauder estimates of Theorem \ref{schauder1}. The existence of $\GK$ then follows from the Ascoli-Arzel\`a's theorem and a diagonal argument. The proof of the existence of a fundamental solution $\GL$ for $\L$ is analogous, and relies on Theorem \ref{ex-propL}.
 
\subsection{Existence of the fundamental solution for the operator $\K$} 
{\sc Proof of Theorem} \ref{main1} \emph{(Existence of the fundamental solution)}. 
We construct a sequence of operators $\left( \K_{n} \right)_{n \in \N}$ satisfying the hypotheses of Theorem \ref{ex-prop}. 
In particular, we need the coefficients $a_n, b_n, r_n$ to be uniformly H\"older continuous and satisfying the assumption \eqref{Hmix}. For this reason, we introduce a cut-off function $\chi_n \in C^{\infty}_{0} (\R^{3})$ such that $0 \le \chi_n (x,y,t) \le 1$, and
\begin{align*}
	\chi_n(x,y,t) = 1 \quad \text{for} \quad x^2 + y^2 \le n^2, \qquad 
	\chi_n(x,y,t) = 0 \quad \text{for} \quad x^2 + y^2 \ge (n+1)^2.
\end{align*}
For every $n \in \N$ we set 
\begin{align*} 
    a_{n}(x,y,t) \, &:= \chi_n(x,y,t) a(x,y,t) + (1- \chi_n(x,y,t)) \lambda, \\ \nonumber
    b_{n}(x,y,t) \, &:= \chi_n(x,y,t) b(x,y,t), \quad \quad
    r_{n}(x,y,t) := \chi_n(x,y,t) r(x,y,t).
\end{align*}
Then we apply Theorem \ref{ex-prop} to the operator $\K_{n}$ for every $n \in \N$. Thus, there exists a sequence of 
equibounded fundamental solutions $(\GK^{\!\!\! n})_{n \in \N}$, in the sense that each of them satisfies \eqref{boundK}.

We define a sequence of open subsets $(\O_{p})_{p \in \N}$ of $\R^6$ such that $\O_{p} \subset \subset \O_{p+1}$ for every $p \in \N$ and $\bigcup_{p=1}^{+ \infty} \O_{p} = \big\{ (x,y,t; \xi, \y,\t) \in \R^{6} \mid  (x,y,t) \ne (\x,\y,\t) \big\}$: 
\begin{equation*}
	\Omega_p :=
	\begin{Bmatrix} 
		(x,y,t; \x, \y, \t) \in \R^6 \mid x^{2} + y^{2} + t^{2} < p^{2}, \x^{2} + \y^{2} + \t^{2} < p^{2}, &  \\ 
		(x-\x)^2 + (y-\y)^2 + (t-\t)^{2} > \left( \tfrac{1}{p} \right)^{2}  &
	\end{Bmatrix} .
\end{equation*} 
We note that the sequence $(\GK^{\!\!\! n})_{n \ge 2}$ is equicontinuous in $\O_{1}$ thanks to Theorem \ref{schauder1}. Moreover, by Theorem \ref{ex-prop} and Theorem \ref{schauder1}, we also have that 
\begin{align*}
	\left(\tfrac{\p \GK^{\!\!\! n}}{\p x} \right)_{n \ge 2}, \quad 
	\left(\tfrac{\p \GK^{\!\!\! n}}{\p \x} \right)_{n \ge 2}, \quad
	\left(\tfrac{\p^{2} \GK^{\!\!\! n}}{\p x^{2}} \right)_{n \ge 2}, \quad 
	\left(\tfrac{\p^{2} \GK^{\!\!\! n}}{\p \x^{2}} \right)_{n \ge 2}, \quad 
	\left(Y \GK^{\!\!\! n} \right)_{n \ge 2}, \quad 
	\text{and} \quad \left(Y^*_{(\x,\y,\t)} \GK^{\!\!\! n} \right)_{n \ge 2} 
\end{align*}
are bounded sequences in $C^{\a}(\O_{1})$. Here $Y$ is the Lie derivative defined in \eqref{lie-diff} and $Y^*_{(\x,\y,\t)}$ is its adjoint, computed with respect to the variable $(\x,\y,\t)$. Thus, there exists a subsequence $(\GK^{\!\!\! 1,m})_{m \in \N}$ that converges uniformly to some function $\G_{1}$ that satisfies \eqref{boundK}. Moreover, $\G_{1} \in C^{2+\a}(\O_{1})$ and, for every $(x_0,y_0,t_0) \in \R^3$ such that $x_0^2 + y_0^2 + t_0^2 <1$ the function $u(x,y,t) := \G_{1}(x,y,t; x_0,y_0,t_0)$ is a classical solution to $\K u = 0$ in the set $\big\{ (x,y,t) \in \R^3 \mid (x,y,t; x_0,y_0,t_0) \in \O_{1} \big\}$, and the function $v(\x,\y,\t) := \G_{1}(x_0,y_0,t_0;\x,\y,\t)$ is a classical solution to $\K^* v = 0$ in the set $\big\{ (\x,\y,\t) \in \R^3 \mid (x_0,y_0,t_0; \x,\y,\t) \in \O_{1} \big\}$.

We next apply the same argument to the sequence $(\GK^{\!\!\! 1,m})_{m \in \N}$ on the set $\O_{2}$, and
obtain a subsequence $(\GK^{\!\!\! 2,m})_{m \in \N}$ that converges in $C^{2+\a}(\O_{2})$ to some function $\G_{2}$, that belongs to $C^{2+\a}(\O_{2})$ and satisfies the bounds \eqref{boundK}, and the following condition. For every $(x_0,y_0,t_0) \in \R^3$ such that $x_0^2 + y_0^2 + t_0^2 <4$ the function $u(x,y,t) := \G_{2}(x,y,t; x_0,y_0,t_0)$ is a classical solution to $\K u = 0$ in the set $\big\{ (x,y,t) \in \R^3 \mid (x,y,t; x_0,y_0,t_0) \in \O_{2} \big\}$, and the function $v(\x,\y,\t) := \G_{2}(x_0,y_0,t_0;\x,\y,\t)$ is a classical solution to $\K^* v = 0$ in the set $\big\{ (\x,\y,\t) \in \R^3 \mid (x_0,y_0,t_0; \x,\y,\t) \in \O_{2} \big\}$. We remark that, since $\G_{2}$ is the limit of a subsequence of $(\GK^{\!\!\! 1,m})_{m \in \N}$, it must coincide with $\G_{1}$ in $\O_{1}$.

We next proceed by induction. Let us assume that the sequence $(\GK^{\!\! q-1,m})_{m \in \N}$ on the set $\O_{q}$ has been defined for some $q \in \N$. We extract from it a subsequence $(\GK^{\!\! q,m})_{m \in \N}$ converging in $C^{2+\a}(\O_{q})$ to some function $\G_{q}$, satisfying \eqref{boundK}. Moreover, $(x,y,t) \mapsto \G_{q}(x,y,t;x_0,y_0,t_0)$ is a classical solution to $\K u = 0$ and $(\x,\y,\t) \mapsto \G_{q}(x_0,y_0,t_0; \x,\y,\t)$ is a classical solution to $\K^* v = 0$. Moreover, it agrees with $\G_{q-1}$ on the set $\O_{q-1}$.

Next, we define a function $\GK$ in the following way: for every $(x,y,t), (\x,\y,\t)  \in \R^{3}$ with $(x,y,t) \ne (\x,\y,\t)$ we choose $q \in \N$ such that $(x,y,t; \x,\y,\t) \in \O_{q}$ and we set $\GK(x,y,t;\x,\y,\t) := \G_{q} (x,y,t;\x,\y,\t)$. This definition is well-posed, since if $(x,y,t) \in \O_{p}$, then $\G_{p}(x,y,t;\x,\y,\t) = \G_{q}(x,y,t;\x,\y,\t)$. 

We next check that $\GK$ has the properties listed in the statement of the Theorem \ref{main1}. As every $\GK^{\!\!\! n}(x,y,t; x_{0}, y_{0}, t_{0})=0$ whenever $t \le t_{0}$, also $\GK(x,y,t; x_{0}, y_{0}, t_{0}) = 0$ for every $t \le t_{0}$. For the same reason, it satisfies \eqref{boundK}. 
Moreover, for every $(x_0,y_0,t_0) \in \R^3$, $(x,y,t) \mapsto \GK(x,y,t; x_0,y_0,t_0) \in L^{1}_{\loc}(\R^{3}) \cap C^{2+\a}_{\loc}(\R^{3} \setminus \{ (x_{0}, y_{0}, t_{0})\})$, and is a classical solution to $\K u = 0$ in $\R^{3} \setminus \{ (x_{0}, y_{0}, t_{0})\}$. Analogously, $(\x,\y,\t) \mapsto \GK(x_0,y_0,t_0;\x,\y,\t) \in L^{1}_{\loc}(\R^{3}) \cap C^{2+\a}_{\loc}(\R^{3} \setminus \{ (x_{0}, y_{0}, t_{0})\})$ and is a classical solution to $\K^* v = 0$ in $\R^{3} \setminus \{ (x_{0}, y_{0}, t_{0})\}$. This proves the point \emph{1.} of the Definition \ref{funsolk} and the point \emph{1.} of Theorem \ref{main1}. 
We remark that points \emph{3.} and \emph{4.} of Theorem \ref{main1} follow immediately from the construction of the fundamental solution $\GK$ and the pointwise convergence.

As far as we are concerned with the reproduction property \emph{2.} of Theorem \ref{main1}, we use the upper bound in \eqref{boundK}, which yields
\begin{equation*}
	\GK^{\!\!\!n} (x,y,t;  \x,\y,\t ) \GK^{\!\!\!n} (\x,\y,\t;  x_{0}, y_{0}, t_{0})  \le 
	C^{+} \, \GK^{\! \! \! \! \l^{+}} (x,y,t; \x,\y,\t ) C^{+} \, \GK^{\! \! \! \! \l^{+}} (\x,\y,\t;  x_{0}, y_{0}, t_{0}),
\end{equation*}
and the reproduction property 
\begin{equation*}
    \int_{\R^2} \GK^{\! \! \! \! \l^{+}} (x,y,t; \x,\y,\t ) 
	\GK^{\! \! \! \! \l^{+}} (\x ,\y ,\t;  x_{0}, y_{0}, t_{0}) d\x \, d \y \, d \t 
	= \GK^{\! \! \! \! \l^{+}} (x,y,t; x_{0}, y_{0}, t_{0}) \, < + \infty,
\end{equation*}
which allows us to use the Lebesgue convergence theorem. Thus the property \emph{2.} of Theorem \ref{main1} holds true.

To proceed with the proof of Theorem \ref{main1} we have to verify that, for every $\phi \in C_{b}(\R^{2})$, the function
	\begin{equation} \label{eq-repr-sol}
		u(x,y,t) \, = \, \int \limits_{\R^{2}} \GK(x,y,t; \x, \y, t_{0}) \, \phi(\x, \y) \, d\x \, d\y
	\end{equation}
is a classical solution to the Cauchy problem
	\begin{equation} \label{eq-cauchy}
		\begin{cases}
			\K u = 0, \qquad &(x,y,t) \in \R^{2} \times ]t_{0}, \infty[; \\
			u(x,y,t_{0}) = \phi(x,y)  &(x,y) \in \R^{2}.
		\end{cases}
	\end{equation}
By the usual standard argument, we differentiate under the integral sign
\begin{equation*}
		\K u(x,y,t) \, = \, \int \limits_{\R^{2}} \K	 \GK(x,y,t; \x, \y, t_{0}) \, \phi(\x, \y) \, d\x \, d\y \, = \, 0.
\end{equation*}
Thus, we are left with the proof that $u(x,y,t) \rightarrow \phi(x_{0},y_{0})$ as $(x,y,t) \rightarrow (x_{0},y_{0},t_{0})$. We first note that
\begin{align*}
	u(x,y,t) - \phi(x_{0}, y_{0}) = &\, \boxed{ \int \limits_{\R^{2}} \GK(x,y,t; \x, \y, t_{0}) 
	\left( \phi(\x, \y) - \phi(x_{0}, y_{0}) \right) d\x \, d\y }_{A} + \\
	&\qquad + \boxed{ \phi(x_{0}, y_{0}) \bigg( \int \limits_{\R^{2}} \GK(x,y,t; \x, \y, t_{0}) d\x \, d\y - 1 \bigg).}_{B}
\end{align*}
The term $B$ plainly vanishes, as $t$ goes to $t_0$, because of the bound \eqref{integralK}. 
The integral $A$ can be handled considering that \eqref{boundK} holds true, and that the expression of $\GK^{\! \! \! \! \l^{+}}$
is \eqref{Glambda}. Specifically, we apply the change of variable 
\begin{equation*}
	\overline x = \tfrac{1}{2 \sqrt{\l^{+} \left(t - t_{0}\right)}} \left( x - \x \right) \quad \quad
	\overline y = \tfrac{\sqrt{3}}{\sqrt{\l^{+} \left(t - t_{0} \right)^{3}}} \left( y - \y + \left(t - t_{0} \right) \tfrac{x + \x}{2} \right),
\end{equation*}
and we obtain the following bound
\begin{equation} \label{eq-expgrowth}
	| A | \, \le \, \frac{1}{\pi}  \int \limits_{\R^{2}} 
	e^{ - \left( \overline{x}^{2} + \overline{y}^{2} \right) } 
	\abs{ \widetilde{\phi} \left( \overline x, \overline y \right) - \phi(x_{0}, y_{0}) }\, d\overline{x} \, d\overline{y},
\end{equation}
where
\begin{equation*}
	\widetilde{\phi} \left( \overline x, \overline y \right) := \phi \left( x - 2 \sqrt{\l^{+} \left(t - t_{0}\right)} \, \overline{x} , \, \,
			y - \sqrt{\tfrac{\l^{+} \left(t - t_{0}\right)^{3}}{3}} \, \overline{y} + \left(t - t_{0} \right) 
			\tfrac{x - \sqrt{\l^{+} \left(t - t_{0}\right)} \, \overline x}{\sqrt{3}} 
	 \right).
\end{equation*}
Note that, for every fixed $(\overline x , \overline y)$, the above expression converges to $\phi(x_{0}, y_{0})$
as $(x,y,t) \rightarrow (x_{0},y_{0},t_{0})$. Moreover $\widetilde{\phi}- \phi(x_0,y_0)$ is bounded as a 
function of $(\overline x, \overline y)$, then the Lebesgue theorem implies that $A$ vanishes as $(x,y,t) \rightarrow (x_{0},y_{0},t_{0})$. Thus $u(x,y,t) \rightarrow \phi(x_{0},y_{0})$ as $(x,y,t) \rightarrow (x_{0},y_{0},t_{0})$, and the proof of the point \emph{1.} of Definition \ref{funsolk} accomplilshed. 

\noindent
The proof that $\GK^{*}(\x, \y, \t; x, y, t) = \GK(x,y,t; \x, \y, \t)$ is the fundamental solution of the adjoint operator 
$\K^{*}$ and satisfies the properties of Theorem \ref{main1} follows accordingly. 
$\hfill \square$

\begin{remark} \label{rem-growht}
The growth condition \eqref{eq-growht} can be used instead of the boundedness assumption on the initial data $\phi$. Indeed, it ensures, alongside with the upper bound \eqref{boundK} for the fundamental solution $\GK$, that the integral \eqref{eq-repr-sol} is convergent for every $(x,y,t) \in \R^{2} \times ]t_{0}, + \infty[$, that it can be differentiated twice with respect to $x$ and once in the direction of the vector field $x \partial_y - \partial_t$, so that $u$ a solution to $\L u = 0$. Moreover, the condition \eqref{eq-growht} and the inequality \eqref{eq-expgrowth} yield that the expression $|A|$ vanishes, as $(x,y,t) \to (x_0, y_0, t_0)$.
%
%
\end{remark}

\subsection{Uniqueness and comparison principle for the operator $\K$}
Now, we recall a technical result, an a priori estimate for the derivatives of the fundamental solution 
$\GK(x,y,t; \x, \y, \t) $ necessary to conclude the proof of Proposition \ref{comp-princ}.
In order to state our result, we introduce for every $R>1$ the set 
\begin{equation}
	\label{set}
	\widetilde Q_{R} := \big\{ R \le \sqrt{\x^2 + \y^2} \le R+1 \big\} \subset \R^{2}.
\end{equation}
\begin{lemma}
	\label{caccioppoliK}
	Let \textbf{(H$_{K}$)} hold, and let $\GK$ be the fundamental solution for the operator $\K$ in the sense of Definition
	\ref{funsolk}. Let $Q_{R}$ be the cylinder defined in \eqref{cylinder} , there exists a constant $C$ such that 
        \begin{align*}
    	 	\frac{\l}{2} \, \int \limits_{t_{0}}^{t} \int_{\widetilde Q_{R}} & \Bigg | \frac{\p \GK}{\p \x}(x,y,t; \x, \y, \t) 
		\Bigg |^{2} d\x d\y d\t  \, \le \\ \nonumber
    	&\le \,  C_{3} \int \limits_{t_{0}}^{t} \int_{\widetilde Q_{R}} \GK(x,y,t; \x, \y, \t)^{2} d\x d\y d\t  \,         
    + \, \frac{1}{2} \int_{\widetilde Q_R} \GK(x,y,t; \x, \y, t_{0})^{2} d\x d\y.
   	 \end{align*}
   	where $C_{3}$ is a positive constant only depending on $\l$, $\Lambda, C_{0}$ and $C_{1}$.
\end{lemma}
The proof of this a priori estimate, also known as Caccioppoli inequality, is presented at the end of this section and 
is based on the representation formula for solutions to the equation $\K u = 0$. For further applications of this technique
see for instance \cite{PP}, \cite{CPP} and \cite{APR}. 

As a first step in the proof of our uniqueness result, we state and prove a comparison principle for the operator $\K$.
\begin{proposition} \label{comp-princ}
Let us consider the operator $\K$ under the assumption \textbf{(H$_{K}$)}. Let $u$ be a classical solution to 
 	\begin{equation} \label{eq-cauchy}
		\begin{cases}
			\K u \ge 0, \qquad &(x,y,t) \in \R^{2} \times ]t_{0}, T]; \\
			u(x,y,t_{0}) \le 0  &(x,y) \in \R^{2}.
		\end{cases}
	\end{equation} 
in the sense of Definition \ref{solution}. 
If moreover
\begin{equation*}
	|u(x,y,t)| \, \le \,M e^{C (x^{2} + y^{2})},
\end{equation*}
for some positive consants $M$ and $C$, then $u \le 0$ in $\R^{2} \times ]t_{0}, T]$.
\end{proposition}
\begin{proof}
We fix a positive constant $\overline t$ such that $\overline t \le 1$, and we prove that, if we choose $\overline t$ small enough, we have $u = 0$ in $\R^2 \times ]t_0, t_0 + \overline t \,]$. We then iterate our argument on the strip $\R^2 \times ]t_0 + \overline t, t_0 + 2 \overline t \,]$, then on $\R^2 \times ]t_0 + 2 \overline t, t_0 + 3 \overline t \,]$. As the choice of $\overline t$ only depends on the operator $\K$ and on the constant $C$ in our assumption $|u(x,y,t)| \, \le \,M e^{C (x^{2} + y^{2})},$ after a finite number of steps we cover the whole set $\R^{2} \times ]t_{0}, T]$. 

Fix $(x,y,s) \in \R^2 \times ]t_0, t_0 + \overline t \,]$ and, denote by $|(x,y)|$ the Euclidean norm of $(x,y)$. For every $R > |(x,y)|$, we let $h_R$ be a $C^\infty(\R^2)$ smooth function, such that $0 \le h_R \le 1$, $h_R (\x,\y) = 1$ whenever $|(\xi, \y)| \le R$, and $h_R (\x,\y) = 0$ for every $(\xi, \y) \in \R^{2}$ with $|(\xi, \y)| \ge R + 1$. We also assume that its first and second order derivatives are bounded uniformly with respect to $R$.

We next recall the Green identity
\begin{equation*}
 v \K u -u \K^{\ast}v = \tfrac{\p}{\p x} \left( v a \tfrac{\p u}{\p x} - u a \tfrac{\p v}{\p x} + b u v \right) + x \tfrac{\p}{\p y} \left( u v \right) - \tfrac{\p}{\p t} \left( u v \right).
\end{equation*}
We then choose a constant $\delta \in ]0, s - t_0[$ and we apply the divergence theorem with $v_{R}(\xi, \eta, \tau) = h_R (\x, \y ) \GK (x, y, s; \xi, \eta, \tau )$, to the cylinder 
\begin{equation*}
 Q_{R,\delta} := \big\{ (\x,\y,\t) \in \R^2 \times ]t_0, t_0 + \overline t \,] \mid |(\xi, \y)| \le R + 2, \tau \le s - \delta \big\}.
\end{equation*}
As $v_{R}, \tfrac{\p v_{R}}{\p x}$, and $\tfrac{\p v_{R}}{\p y}$ vanish at the lateral part of the boundary of $Q_{R,\delta}$, we find
\begin{equation} \label{eq-green}
  \int_{Q_{R,\delta}} \left( v_{R} \K u -u \K^{\ast}v_{R} \right) (\x,\y,\t) d\x \, d\y\, d\t 
  = - \int_{\R^2} \left( u v_{R} \right) (\x,\y, s - \delta) d\x \, d\y + \int_{\R^2} \left( u v_{R} \right) (\x,\y, 0) d\x \, d\y.
\end{equation}

Because of the properties of the fundamental solution we have
\begin{equation*}
 u(x,y,s) = \lim_{\delta \to 0} \int_{\R^2}  \GK (x, y, s, \xi, \eta, s - \delta ) h_R (\x, \y ) u (\x,\y, s - \delta) d\x \, d\y.
\end{equation*}
Moreover, by our assumption, we have $v_{R} \K u \ge 0$ in $\R^2 \times ]t_0,t_1[$, and $(uv_{R})(\cdot, \cdot, 0) \le 0$. Hence \eqref{eq-green} gives
\begin{align*}
 u(x,y,s) = & \int_{\R^2} \left( u v_{R} \right) (\x,\y, 0) d\x \, d\y - 
 \int_{\R^2 \times ]t_0,s[} \left( v_{R} \K u - u \K^{\ast}v_{R} \right) (\x,\y,\t) d\x \, d\y\, d\t \\
 \le & \int_{\R^2 \times ]t_0,s[} u (\x,\y,\t) \K^{\ast}v_{R}  (\x,\y,\t) d\x \, d\y\, d\t.
\end{align*}

We are left with the proof that the right hand side of the above inequality vanishes as $R \to + \infty$.
From now on, we only sketch the proof since it suffices to proceed as in the proof of Theorem 1.6 of \cite{PDF}. 
Since $\K^{\ast} \GK (x,y,s;\x,\y,\t) = 0$, we deduce
\begin{align*}
 u(x,y,s)  \le & 2 \Lambda \int \limits_{t_0}^{s} \bigg(\int_{\widetilde Q_{R}} |u (\x,\y,\t)|  
 \left| \frac{\p \GK}{\p \x} (x,y,s; \x,\y,\t) \right| \left| \frac{\p h_R}{\p \x} (\x,\y,\t)\right|  d\x \, d\y \bigg) d\t + \\ 
 &+ \int\limits_{t_0}^{s} \bigg(\int_{\widetilde Q_{R}} |u (\x,\y,\t)| \left| \GK (x,y,s; \x,\y,\t) \right| |\K^{\ast} h_R (\x,\y,\t) |  d\x \, d\y \bigg) d\t,
\end{align*}
where $\widetilde Q_{R}$ is defined in \eqref{set}. 
We recall that first and second order derivatives of the function $h_{R}$ are bounded because of its definition, 
more precisely we have that
\begin{align}
	\label{derivate}
	\left| Y(h_{R}(x,y,t) \right| \, &\le \, C_{0} R, \qquad \left| \frac{\p h_{R}(x,y,t)}{\p x} \right| \, \le \, C_{1},
	 \quad  \text{and}\\ \nonumber
	 \left| K^{*} h_{R}(x,y,t) \right| \, &\le \, \Lambda \left( 2 + C_{1} \right) + C_{0} R =: C_{2}(1 + R),
\end{align} 
where $C_{0}, C_{1}$ and $C_{2}$ are positive constants. Thus, by applying the H\"older inequality and the estimates 
\eqref{derivate}, we get the following inequality
\begin{align*}
 u(x,y,s) \le 2 \Lambda C_{1} & \bigg( \int\limits_{t_0}^{s} \int_{\widetilde Q_{R}} |u (\x,\y,\t)|^{2}  
 				d\x d\y d\t \bigg)^{\frac12}
 			 \bigg( \int\limits_{t_0}^{s} \int_{\widetilde Q_{R}} \left|
			 \frac{\p \GK}{\p \x} (x,y,s; \x,\y,\t) \right|^2 d\x d\y d\t \bigg)^{\frac12} \, + \\ 
 		C_{2}(1 + R) & \bigg( \int\limits_{t_0}^{s} \int_{\widetilde Q_{R}} |u (\x,\y,\t)|^{2} d\x d\y d\t \bigg)^{\frac12}
		     \bigg( \int\limits_{t_0}^{s} \int_{\widetilde Q_{R}} \GK (x,y,s; \x,\y,\t)^{2}  d\x d\y d\t \bigg)^{\frac12}.
\end{align*} 
Then, Lemma \ref{caccioppoliK} yields 
\begin{align*}
 u(x,y,s) \le & \left(2 \Lambda C_{1} C_3 + C_2 (1+R) \right)  
 \bigg( \int\limits_{t_0}^{s} \int_{\widetilde Q_{R}} |u (\x,\y,\t)|^{2} d\x d\y d\t \bigg)^{\frac12} \cdot \\
        & \bigg( \int\limits_{t_0}^{s} \int_{\widetilde Q_{R}} \GK (x,y,s; \x,\y,\t)^{2}  d\x d\y d\t +
        \int_{\widetilde Q_R} \GK(x,y,s; \x, \y, t_{0})^{2} d\x d\y \bigg)^{\frac12}.
\end{align*} 
By our assumption $|u(\x,\y,\t)| \, \le \,M e^{C (\x^{2} + \y^{2})}$, we have that 
\begin{equation*} 
	\bigg( \int\limits_{t_0}^{s} \int_{\widetilde Q_{R}} |u (\x,\y,\t)|^{2} d\x d\y d\t \bigg)^{\frac12} \le 
	2 \sqrt{\pi \overline t R} \, M e^{C (R+1)^2}.
\end{equation*}
Moreover, the Corollary \ref{corollary} gives
\begin{equation*} 
	\bigg( \int\limits_{t_0}^{s} \int_{\widetilde Q_{R}} \GK (x,y,s; \x,\y,\t)^{2}  d\x d\y d\t +
        \int_{\widetilde Q_R} \GK(x,y,s; \x, \y, t_{0})^{2} d\x d\y \bigg)^{\frac12} \le 
         2 \sqrt{\pi (1 + \overline t) R} \overline C e^{- \overline C \frac{(R-1)^2}{\overline t}},
\end{equation*}
provided that $R-1$ is greater than the constant $R_0$ appearing in its statement. 
Finally, recalling that $0 < \overline t \le 1$, we conclude that
\begin{align*}
 u(x,y,s) \le & C_{4} (1+R)^2 e^{C (R+1)^2} e^{- \overline C \frac{(R-1)^2}{\overline t}},
\end{align*} 
where $C_4$ is a positive constant depending on the operator $\K$. In order to conclude our proof, it sufficies to choose $\overline t < \frac{\overline C}{C}$. The concludion the follows by letting $R \to + \infty$. Hence $u(x, y, s) \le 0$. 
The thesis follows by repeating the previous argument finitely many times, as the choice of $\overline t$ does not depend on $(x,y,s)$. 
\end{proof}

\bigskip
\noindent
{\sc Proof of Theorem} \ref{uniquenessK}. This uniqueness result plainly follows from Proposition \ref{comp-princ} firstly applied to 
$u= u_{1} - u_{2}$, and then to $u=u_{2} - u_{1}$.
$\hfill \square$

\bigskip
\noindent
{\sc Proof of Theorem} \ref{main1} \emph{(Uniqueness of the fundamental solution)}. Suppose that $\G_{1}$ and $\G_{2}$ are two fundamental  
solutions for the operator $\K$. For every $\phi \in C^{\infty}_{b}(\R^{2})$ we define 
\begin{equation*} 
		u_{1}(x,y,t) = \int \limits_{\R^{2}} \G_{1}(x,y,t; \x, \y, t_{0}) \, \phi(\x, \y) \, d\x \, d\y, \quad
		u_{2}(x,y,t) = \int \limits_{\R^{2}} \G_{2}(x,y,t; \x, \y, t_{0}) \, \phi(\x, \y) \, d\x \, d\y
	\end{equation*}
and we note that are bounded classical solutions to the same Cauchy problem \eqref{eq-cauchy}. Then 
$u_{1} = u_{2}$ by Theorem \ref{uniquenessK}. Since $\phi$ is arbitrarily chosen we have that $\G_{1} = \G_{2}$.
$\hfill \square$ 

\bigskip
\noindent
{\sc Proof of Lemma} \ref{caccioppoliK}. For every $R > 0$, let us consider the following cylinder
     \begin{equation}
		\label{cylinder}
		 Q_{R} := \big\{ (\x,\y,\t) \in \R^2 \times ]t_0, t \,] \mid |(\xi, \y)| \le R + 2, \tau \le t \big\},
    \end{equation}
    which is a slight modification of the cylinder $Q_{R,\delta}$ previously introduced in the proof of Proposition \ref{comp-princ}.
     Let us consider the fundamental solution $\GK$ associated with the operator $\K$. By definition, $\GK$ satisfies the equation 
     $\K u = 0$. Thus, by multiplying the equation by a certain test function $\phi(x,y,t) \in C^{\infty}_{0}(\R^{3})$, integrating on the 
     cylinder $Q_{R}$ and then proceeding by parts, we get the following equality
     \begin{align*}
    	0 = \, &-  \int_{\Q_R} \Big \langle a  \frac{\p \GK}{\p \x}(x,y,t; \x, \y, \t), \frac{\p  \phi}{\p \x} \Big \rangle 
			\, d\x \, d\y \, d\t \, + \, \int_{\Q_R} b \, \frac{\p \GK}{\p \x}(x,y,t;\x, \y, \t)  \, \phi \, d\x \, d\y \, d\t \, + \\
	   &+ \, \int_{\Q_R} Y \GK(x,y,t;\x, \y, \t) \phi \, d\x \, d\y \, d\t \, - \, \int_{\Q_R} r \GK(x,y,t;\x, \y, \t)  \, \phi \, d\x \, d\y \, d\t .
    \end{align*}
    This equality is also known as \textit{weak formulation} of the equation $\K u = 0$, and 
    for the sake of clarity from now on we set $\GK = \GK(x,y,t; \x, \y, \t)$. 
    In particular, as a test function we can consider
    $\phi(\x,\y,\t) := \left[ h_{R+1}(\x, \y) - h_{R-1}(\x, \y) \right]^{2} \GK$, where 
    $h_R$ is the same smooth function introduced in the proof of Proposition \ref{comp-princ}, we get 
    \begin{align*}
    	0 \le \phi \le 1, \quad \phi = \begin{cases}
							0 \quad &\text{for} \quad  |(\xi, \y)| \le R-1, \\ 
							(1 - h_{R-1})^{2}\, \GK \quad &\text{for} \quad  R -1 < |(\xi, \y)| \le R \\
						        \GK \quad &\text{for} \quad  R < |(\xi, \y)| \le R + 1 \\
        				        			h_{R+1}^{2} \, \GK \quad &\text{for} \quad R+1 < |(\xi, \y)| < R + 2 \\
							0 \quad &\text{for} \quad |(\xi, \y)| \ge R +2.
						\end{cases}
    \end{align*}
    Since $\p_{\x} \phi = \left( h_{R+1} - h_{R} \right) \p_{\x} \GK + 
    2 \left( h_{R+1} - h_{R} \right) \GK \p_{\x}\left( h_{R+1} - h_{R}\right)$,  
    assumption \textbf{(H$_{K}$)} holds true and the first and second order derivatives of the function $h_{R}$ are bounded as 
    in \eqref{derivate}, we get the following inequality
    \begin{align}
        \label{estimate}
        \l \, \int \limits_{t_{0}}^{t}  \int \limits_{\widetilde Q_{R}} \left| \frac{\p \GK}{\p \x} \right|^2  \, 
        &\le \, \boxed{2 \l \int_{\Q_R} \left| \GK \, \left( h_{R+1} - h_{R-1} \right) \, \frac{\p \GK}{\p \x} \, 
        \frac{\p}{\p \x} \left( h_{R+1} - h_{R-1} \right) \right| \,}_{A} 
        \, + \\ \nonumber
        &+ \, \boxed{ \frac{1}{2} \int_{\Q_R} Y\left( \GK^{2} \right) \left(h_{R+1}(\x, \y) - h_{R-1}(\x, \y) \right)^{2}}_B \, 
        \\ \nonumber
        &+ \, \boxed{\Lambda \int_{\Q_R} \left| \GK \, \frac{\p \GK}{\p \x} \, \left( h_{R+1} - h_{R-1} \right)^{2} \right| }_{C} \, 
         +  \, \Lambda \int_{\Q_R} \GK^{2} \, \left( h_{R+1} - h_{R-1} \right)^{2},
    \end{align}
    where the set $\widetilde Q_{R}$ has previously been defined in \eqref{set}.
Now, we can estimate terms A and C by Young's inequality. 
As far as we are concerned with term B, we begin considering the following identity:
\begin{equation*}
	\left( h_{R+1} - h_{R} \right)^{2} Y(\GK^{2}) \, = \, Y \left(\GK^{2} \left( h_{R+1} - h_{R} \right)^{2}\right) \, 
	- \, \GK^{2} \, Y\left(\left( h_{R+1} - h_{R} \right)^{2}\right).
\end{equation*} 
Thus, we can rewrite term B as the sum of two terms, and by applying the divergence theorem to B$_{1}$
($\GK^{2} \left( h_{R+1} - h_{R} \right)^{2}$ is null on the lateral boundary of $\widetilde Q_{R}$), we get
\begin{align*}
	&\boxed{ \frac{1}{2} \int_{\Q_R} Y\left( \GK^{2} \right) \left(h_{R+1}(\x, \y) - h_{R}(\x, \y) \right)^{2} }_B \, \le  \\	
	&\le \boxed{\frac{1}{2} \int_{\Q_R} Y \left(\GK^{2} \left( h_{R+1} - h_{R} \right)^{2}\right) \, }_{B_{1}}
	    + \, \boxed{\int_{\Q_R} \GK^{2} \, Y\left(\left( h_{R+1} - h_{R} \right)^{2}\right)}_{B_{2}} \\ 
         &\le \, \frac{1}{2} \int_{\widetilde Q_R} \GK (x,y,t; \x, \y, t_{0})^{2} d\x d\y \, + \,
         	2 C_{0} \int \limits_{t_{0}}^{t}  \int_{\widetilde Q_{R}} \GK^{2} d\x d\y d\t.
\end{align*}
By choosing $\e = \frac{\lambda}{2 \left( 4 \l C_{1} - \Lambda \right)}$ we get
\begin{align*}
    \frac{\l}{2} \, \int \limits_{t_{0}}^{t} & \int_{\widetilde Q_{R}} \Bigg | \frac{\p \GK}{\p \x} \Bigg |^{2} d\x d\y d\t  \,
    \le \, C_{3} \int \limits_{t_{0}}^{t} \int_{\widetilde Q_{R}} \GK ^{2} d\x d\y d\t  \,         
    + \, \frac{1}{2} \int_{\widetilde Q_R} \GK^{2}(x,y,t; \x, \y, t_{0}) d\x d\y,
\end{align*}
where $C_{3}= C_{3}(\l, \Lambda, C_{0}, C_{1})$ is a positive constant.
$\hfill \square$

\subsection{Existence of the fundamental solution for the operator $\L$}
{\sc Proof of Theorem} \ref{main2} \emph{(Existence of the fundamental solution)}. The proof of this theorem is analogous to the proof of Theorem \ref{main1}. In this case, 
we construct a sequence of operators $\left( \L_{n} \right)_{n \in \N}$ satisfying the assumptions of Theorem \ref{ex-propL}. 
In particular, we need the coefficients $a_n, b_n$ to be smooth and satisfying a suitable version of the condition 
\eqref{HmixL}. For this reason, we introduce a non-negative function $\r \in C^{\infty}_{0} (\R^{3})$ such that
\begin{align*}
	\int_{\R^{3}} \r = 1, \qquad 	B_{0} := \text{\rm supp} \, 
	\r \subset \left\{ (x,y,t) \in \R^3 \mid x^2 + y^2 + t^2 < \tfrac14 \right\}, 
\end{align*}
and then proceed with a standard mollifying procedure. 
In order to take into consideration the fact that the domain of the coefficients $a$ and $b$ is $\R^+ \times \R^2$, for every 
$(x,y,t) \in \R^+ \times \R^2$ and for every $n \in \N$ we set
\begin{align*} 
    a_{n}(x,y,t) \, &:= \, \int \limits_{B_{0}} a \left(x - \tfrac{x \x}{n}, y - \tfrac{\y}{n} , 
    t - \tfrac{\t}{n}\right) \, \r(\x, \y, \t) \, d\x \, d\y \, d\t, \\ \nonumber
    b_{n}(x,y,t) \, &:= \, \int \limits_{B_{0} } b \left(x - \tfrac{x \x}{n}, y - \tfrac{\y}{n} , 
    t - \tfrac{\t}{n}\right) \, \r(\x, \y, \t) \, d\x \, d\y \, d\t.
\end{align*}
Note that $ \left(x - \tfrac{x \x}{n}, y - \tfrac{\y}{n} , t - \tfrac{\t}{n}\right) \in B(x,y,t)$ for every 
$(\x, \y, \t) \in B_{0}$ and for every $n \in \N$, where
\begin{equation*}
	B(x,y,t) := \left[\tfrac{1}{2} x, \tfrac32 x \right] \times  \left[y - \tfrac{1}{2}, y + \tfrac{1}{2} \right] \times 
	\left[t - \tfrac{1}{2}, t + \tfrac{1}{2} \right].
\end{equation*}
Then for every $n \in \N$ the coefficients $a_{n}$ are smooth and satisfy the following version of 
\textbf{(H$_{L}$)}
\begin{align*}
	\lvert a_{n}(x,y,t) \rvert \le \sup \limits_{B(x,y,t)} |a| \le \, \Lambda, \quad 
	\Big \lvert \tfrac{\p a_{n} }{\p x} (x,y,t) \Big \rvert \le \sup \limits_{B(x,y,t)} \Big \lvert \tfrac{\p a }{\p x} \Big \rvert 
			      \le \, \tfrac{2\Lambda}{x},  \quad 
	a_{n}(x,y,t) \ge \, \inf \limits_{B(x,y,t)} a
			  \ge \, \lambda.
\end{align*}
The same statement holds true for the coefficients $b_{n}$, with $n \in \N$. 
Then we apply Theorem \ref{ex-propL} to the operator $\L_{n}$ for every $n \in \N$. Thus, there exists a sequence of 
equibounded fundamental solutions $(\GL^{\!\!\! n})_{n \in \N}$, in the sense that each of them satisfies \eqref{boundL}.

Then we apply the same diagonal argument as in the proof of Theorem \ref{main1}, but with a different choice for 
the open sets $(\O_{p})_{p \in \N}$ of $(\R^+ \times \R^{2})^{2}$. Indeed, we define
\begin{equation*}
	\Omega_p := \left\{  \begin{matrix}
				(x,y,t; \x, \y, \t) \in (\R^+ \times \R^{2})^{2} \mid x^{2} + y^{2} + t^{2} \le p^{2}, \quad \x^{2} + \y^{2} + \t^{2} \le p^{2}  \\
		 		(x - \x)^{2} + (y - \y)^{2} + (t - \t)^{2} \ge \frac{1}{2p}, \quad
				x > \frac{1}{p}, \quad \x > \frac{1}{p}
				\end{matrix} \right\},
\end{equation*} 
such that $\bigcup_{p=1}^{+ \infty} \O_{p} = \big\{ (x,y,t; \xi, \y,\t) \in (\R^+ \times \R^{2})^{2} \mid  (x,y,t) \ne (\x,\y,\t) \big\}$ and
$\O_{p} \subset \subset \O_{p+1}$ for every $p \in \N$. Thus, we define a function $\GL$ in the following way: 
for every $(x,y,t), (\x, \y, \t) \in \R^+ \times \R^{2}$ with $(x,y,t) \ne (\x, \y, \t)$ we choose $q \in \N$ such that $(x,y,t; \x, \y, \t) \in \O_{q}$
and we set $\GL(x,y,t;\x,\y,\t) := \G_{q} (x,y,t;\x,\y,\t)$. This definition is well-posed, since if $(x,y,t) \in \O_{p}$, then $\G_{p}(x,y,t;\x,\y,\t) = \G_{q}(x,y,t;\x,\y,\t)$. 

We next check that $\GL$ has the properties listed in the statement of the Theorem \ref{main2}. As every $\GL^{\!\!\! n}(x,y,t; x_{0}, y_{0}, t_{0})=0$ whenever $t \le t_{0}$ or $y \ge y_{0}$, also $\GL(x,y,t; x_{0}, y_{0}, t_{0}) = 0$ whenever $t \le t_{0}$ or $y \ge y_{0}$. 
For the same reason, it satisfies \eqref{boundL}. 
Moreover, for every $(x_0,y_0,t_0) \in \R^+ \times \R^{2}$, $(x,y,t) \mapsto \GL(x,y,t; x_0,y_0,t_0) \in L^{1}_{\loc}(\R^+ \times \R^{2}) 
\cap C^{2+\a}_{\loc}(\R^+ \times \R^{2} \setminus \{ (x_{0}, y_{0}, t_{0})\})$, and is a classical solution to $\L u = 0$ in 
$\R^+ \times \R^{2} \setminus \{ (x_{0}, y_{0}, t_{0})\}$. Analogously, $(\x,\y,\t) \mapsto \GL(x_0,y_0,t_0;\x,\y,\t) \in L^{1}_{\loc}(\R^+ \times \R^{2}) \cap C^{2+\a}_{\loc}(\R^{3} \setminus \{ (x_{0}, y_{0}, t_{0})\})$ and is a classical solution to $\L^* v = 0$ in $\R^+ \times \R^{2} \setminus \{ (x_{0}, y_{0}, t_{0})\}$. This proves the point \emph{1.} of the Definition \ref{funsoll} and the point \emph{1.} of Theorem \ref{main2}. 
We remark that points \emph{3.} and \emph{4.} of Theorem \ref{main2} follow immediately from the construction of the fundamental solution $\GL$ 
and the pointwise convergence.
As far as we are concerned with the reproduction property \emph{2.} of Theorem \ref{main2}, 
we proceed as in the proof of Theorem \ref{main1} thanks to Corollary \ref{corollaryL}.

To proceed with the proof of Theorem \ref{main2} we have to verify that for every $\phi \in C_{b}(\R^{2})$ the function
		\begin{equation*}
			u(x,y,t) \, = \, \int \limits_{\R^{2}} \GL(x,y,t; \x, \y, t_{0}) \, \phi(\x, \y) \, d\x \, d\y
		\end{equation*}
is a classical solution to the Cauchy problem
		\begin{equation*}
			\begin{cases}
				\L u (x,y,t) = 0, \qquad &(x,y,t) \in \R^{+} \times \R \times \R^{+}; \\
				u(x,y,t_{0}) = \phi(x,y)  &(x,y) \in \R^{+} \times \R.
			\end{cases}
		\end{equation*}
By a very standard argument we differentiate under the integral sign and we find
\begin{equation*}
		\L u(x,y,t) \, = \, \int \limits_{\R^{+} \times \R} \L \GL(x,y,t; \x, \y, t_{0}) \, \phi(\x, \y) \, d\x \, d\y \, = \, 0.
\end{equation*}
Thus, to conclude the proof we have to verify that for any $(x_{0}, y_{0}) \in \R^{+} \times \R$ we have 
\begin{equation}
	\label{convergenza}
	\lim \limits_{(x,y,t) \rightarrow (x_{0}, y_{0}, 0)} u(x,y,t) = \phi (x_{0}, y_{0}) .
\end{equation} 
The proof of this fact is based on the use of ``barriers'', and on Theorems 6.1 and 6.3 of \cite{M}. 
The following argument relies on the fact that the operator $\L$ behaves as the operator $\K$ in every compact set of $\R^{+} \times \R \times \R$.
Let us consider the sequence of functions
\begin{equation*}
	u_{n}(x,y,t) \, = \, \int \limits_{\R^{2}} \GL^{\! \!n}(x,y,t; \x, \y, t_{0}) \, \phi(\x, \y) \, d\x \, d\y
\end{equation*}
and note that $u(x,y,t) = \lim \limits_{n \rightarrow \infty} u_{n}(x,y,t)$. Since $\GL^{\! \!n}$ is the fundamental solution of $\L_{n}$, we have that 
\begin{equation}
	\label{convergenzan}
	\lim \limits_{(x,y,t) \rightarrow (x_{0}, y_{0}, 0)} u_{n}(x,y,t) = \phi (x_{0}, y_{0})  
	\quad \text{for every} \, n \in \N.
\end{equation} 
Let us introduce the cylinder
\begin{equation*}
	Q:= \left]\tfrac{1}{2}x_{0} , \tfrac{3}{2}x_{0} \right[ \times \left] y_{0} - 1 , y_{0} + 1\right[ \times ] 0, T[ 
\end{equation*} 
centered at $(x_{0}, y_{0}, 0)$. As the sequence $\{ u_{n} \}_{n \in \N}$ is uniformly bounded, it is possible 
to construct two barrier functions $u^{+}$ and $u^{-}$, a super and a sub solution respectively, such that 
\begin{equation*}
	u^{-} (x,y,t) \le u_{n}(x,y,t) \le u^{+} (x,y,t) \quad \text{for every} \, (x,y,t) \in Q
\end{equation*} 
and for every $n \in \N$, and such that 
\begin{equation*}
	\lim \limits_{(x,y,t) \rightarrow (x_{0}, y_{0}, 0)} u^{-}(x,y,t) = \phi (x_{0}, y_{0}) , \qquad 
	\lim \limits_{(x,y,t) \rightarrow (x_{0}, y_{0}, 0)} u^{+}(x,y,t) = \phi (x_{0}, y_{0}).
\end{equation*} 
The claim \eqref{convergenza} directly follows. 
$\hfill \square$

\begin{remark} \label{rem-growht-L} 
The linear growth of the initial condition in the Cauchy problem \eqref{PDE1} is allowed in the formula \eqref{rep2}. Indeed, the 
Corollary \ref{corollaryL} holds for the operator $\L$ satisfying the assumption {\rm \textbf{(H$_{L}$)}}, and it is known that the Geman-Yor process \eqref{processY} has finite first order moments. 
\end{remark}

\subsection{Uniqueness and comparison principle for the operator $\L$}
Following the steps of the proof of Theorem \ref{uniquenessK} for the uniqueness of the solution for the Cauchy problem 
associated to the operator $\K$, we need to prove an intermediate result for the operator $\L$, also known as comparison principle. 
In particular, we can apply the general result due to Aronson and Besala proved in \cite{AronsonBesala}, 
that in the case of the operator  $\L$ reads as follows.

\medskip
\noindent
{\sc Theorem B} \textit{(Aronson - Besala). 
	Let us consider for $T>0$ the open set $\O = \R^{+} \times \R \times ]0,T]$ and let $\L$ be the differential operator 
	defined in \eqref{Ldiv} under the assumption \textbf{(H$_{L}$)}. 
	If $u$ is a classical solution of $\L u \le 0$ in $\O$
	such that 
	\begin{equation*}
		u(x,y,0) \ge 0, \quad \text{for} \, (x,y) \in \R^{+} \times \R
		\quad \text{and} \quad u(0,y,t) \ge 0 \quad \text{for} \, (y,t) \in \R \times  ]0,T]
	\end{equation*}
	and for some positive constant $M$ and $k$
	\begin{equation*}
		u(x,y,t) \ge - M \, \exp \left\{ k \, \log \left( x^{2} + y^{2} + 1 \right) + 1 \right\}^{2}
	\end{equation*}
	in $\O$, then $u(x,y,t) \ge 0$ in $\overline \O$.
}

\medskip
\noindent
We remark that the above results would be enough to ensure the uniqueness of the solution for the Cauchy problem 
associated to the operator $\L$ in the form \eqref{Ldiv}. Nevertheless, when considering the operator $\L$ with locally H\"older continuous coefficients
and satisfying the assumption \textbf{(H$_{L}$)} as in our case, we can improve the previous result by requiring the 
solution $u$ to have a positive sign only on the boundary related to the initial data
\begin{equation*}
	u(x,y,0) \ge 0, \quad \text{for} \, (x,y) \in \R^{+} \times \R,
\end{equation*}
and getting rid of the sign assumption on the part of the boundary $\left\{ 0 \right\} \times \R \times ]0,T]$.
\begin{theorem}
	\label{APtheoremB}
	Let us consider for $T>0$ the open set $\O = \R^{+} \times \R \times ]0,T]$ and let $\L$ be the differential operator 
	defined in \eqref{Ldiv} under the assumption \textbf{(H$_{L}$)}. 
	If $u$ is a classical solution of $\L u \le 0$ in $\O$
	such that 
	\begin{equation} \label{sign}
		u(x,y,0) \ge 0, \quad \text{for} \, (x,y) \in \R^{+} \times \R,
	\end{equation}
	and for some positive constant $M$ and $k$
	\begin{equation} \label{lower}
		u(x,y,t) \ge  - M \exp\left(C (\log(x^{2} + y^{2} + 1) - \log(x)) + 1\right)^2,
	\end{equation}
	for every $(x,y,t) \in \R^{+} \times \R \times ]0, T]$. Then $u \ge 0$ in $\R^{+} \times \R \times [0, T]$.
\end{theorem}
\bigskip
\noindent
{\sc Proof} 
Let us consider, for a given $\b >0$, the auxiliary function 
	\begin{align*}
		v(x,y,t) = \exp \left( 2 e^{\b t} \, C (\log(x^{2} + y^{2} + 1) - \log(x)) + 1\right)^2.
	\end{align*}
	It is easily verified that if $t \in ]0, 1/\b]$ we have
	\begin{align*}
		\L v(x,y,t) \le e^{\b t} v \left(C (\log(x^{2} + y^{2} + 1) - \log(x)) + 1\right)^2 \left( E - 2\b \right),
	\end{align*}
	where $E$ is a positive constant only depending on the constants $C$ and $\l, \Lambda$ appearing in \textbf{(H$_{L}$)}. Thus, 
	if we set $\b = E$ it follows that $\L v < 0$ in $\R^{+} \times \R \times ]0, 1/\b]$.

In the following, we let $\widetilde \psi(x,y):= \log(x^{2} + y^{2} + 1) - \log(x)$ and we note that, 
for every $K > \log(2)$, we have
\begin{equation*}
 \left\{ (x,y) \in \R^+ \times \R \mid \widetilde \psi(x,y) < K \right\} = \left\{ (x,y) \in \R^2 \mid (x- x_K)^2 + y^2 < r_K^2 \right\},
\end{equation*}
where $x_K := \tfrac{e^K}{2}$ and $r_K = \sqrt{x_K^2-1}$. If $\b >0$ is as above, we consider, for arbitrary $K > \log(2)$ and $M >0$, the function
	\begin{align*}
		w(x,y,t) = u(x,y,t) +  M e^{- \left( C K + 1 \right)^{2}}  v(x,y,t), 
	\end{align*}

It is clear that $\L w <0$ in $\R^{+} \times \R \times ]0, 1/\beta]$, and that $w(x,y,0) \ge 0$ for $(x,y) \in \R \times \R^{+}$, by \eqref{sign}. Moreover, because of \eqref{lower}, we have that
	\begin{align*}
		w(x,y,t) \ge 0  \qquad &\text{for} \, (x,y,t) \in 
		\big\{ \R^{+} \times \R \times [0,T] \mid \widetilde \psi(x,y) = K \big\}.
	\end{align*}
From the weak minimum principle it follows that $w(x,y,t) \ge 0$ for every $(x,y,t) \in \R^{+} \times \R \times [0,1/\beta]$ such that $\widetilde \psi(x,y) \le K$.

Now, if $(x, y, t)$ is any point in $\R^{+} \times \R \times ]0, 1/\b]$, we choose $K$ such that $\widetilde \psi(x,y) \le K$, and by the above argument it follows that $w(x,y,t) \ge 0$. The case $t> 1/\b$ straightly follows by repeating the above argument.
$\hfill \square$

\bigskip
\noindent
{\sc Proof of Theorem} \ref{uniquenessL}. This uniqueness result plainly follows from Proposition \ref{APtheoremB} firstly applied to 
	$u= u_{1} - u_{2}$, and then to $u=u_{2} - u_{1}$.
$\hfill \square$

\bigskip
\noindent
{\sc Proof of Theorem} \ref{main2} \emph{(Uniqueness of the fundamental solution)}. Suppose that $\G_{1}$ and $\G_{2}$ are two fundamental  
solutions for the operator $\L$. For every $\phi \in C^{\infty}_{b}(\R^{2})$ we define 
\begin{equation*} 
		u_{1}(x,y,t) = \int \limits_{\R^{2}} \G_{1}(x,y,t; \x, \y, t_{0}) \, \phi(\x, \y) \, d\x \, d\y, \quad
		u_{2}(x,y,t) = \int \limits_{\R^{2}} \G_{2}(x,y,t; \x, \y, t_{0}) \, \phi(\x, \y) \, d\x \, d\y
\end{equation*}
and we note that are bounded classical solutions to the same Cauchy problem. Then 
$u_{1} = u_{2}$ by Theorem \ref{uniquenessL}. Since $\phi$ is arbitrarily chosen we have that $\G_{1} = \G_{2}$.
$\hfill \square$

\bibstyle{siam}

\end{document}